# Maximum pseudolikelihood estimator for exponential family models of marked Gibbs point processes

**Jean-Michel Billiot, Jean-François Coeurjolly and Rémy Drouilhet**

*SAGAG Team, Dept. of Statistics, LJK, Grenoble University, France*
*e-mail:* Jean-Michel.Billiot@upmf-grenoble.fr*;*
Jean-Francois.Coeurjolly@upmf-grenoble.fr*;* Remy.Drouilhet@upmf-grenoble.fr

**Abstract:** This paper is devoted to the estimation of a vector $\boldsymbol{\theta}$ parametrizing an energy function of a Gibbs point process, via the maximum pseudolikelihood method. Strong consistency and asymptotic normality results of this estimator depending on a single realization are presented. In the framework of exponential family models, sufficient conditions are expressed in terms of the local energy function and are verified on a wide variety of examples.

## 1. Introduction

The class of Gibbs point processes is interesting because it allows us to introduce and study interactions between points through the modelling of an associated energy function. Historical aspects of the mathematical theory are covered briefly in Kallenberg (1983). When the energy function is parametrized, one among many other methods of estimation, is the maximization of the pseudolikelihood. Baddeley and Turner (2000) dealt with some practical aspects of such a parametric method and gave a survey on asymptotic results. They noticed that some classical examples (such as the area-interaction model, the Multi-Strauss model, the 2-type Strauss model) do not satisfy the assumptions of the existing asymptotic normality results. This paper aims at filling this gap.

Many proposals tried to estimate the energy function from the available point pattern data generated by some marked Gibbs point processes. If the energy belongs to a parametric family model, the most well-known methodology is the use of the likelihood function, see *e.g.* Møller and Waggepetersen (2003) and the references therein. The main drawback of this approach is that the likelihood function contains an unknown scaling factor whose value depends on the parameters and which is difficult to calculate. An alternative approach relies on the use of the pseudolikelihood. This idea originated from Besag (1974) in the study of lattice processes. Besag et al. (1982) further considered this method for pairwise interaction point processes, while Jensen and Møller (1991) generalized it to the general class of marked Gibbs point processes. A general review of the problem of statistical inference on spatial point processes including the Takacs-Fiksel method (a parametric method based on a characteristic property of marked Gibbs point processes using Palm measure) and non-parametric methods can be found in the recent monograph of Møller and Waggepetersen (2003).

In order to underline our theoretical contributions, let us present the different papers discussing asymptotic properties of the maximum pseudolikelihood estimator. The first work was done by Jensen and Møller (1991). The authors obtained consistency for exponential family models of marked point processes. They mainly considered inhibition and hard-core models. They notably applied their results on the marked Strauss process. Jensen and Künsch (1994) and Mase (1995) focused on specific models with two parameters -the chemical potential and the inverse temperature- which can be viewed as particular exponential family models. Jensen and Künsch (1994) obtained an asymptotic normality result by first assuming the inhibition or hard-core property and then the finite range property. Mase (1995) established consistency for the class of superstable and lower regular potentials introduced by Ruelle (1970). Mase (2000) extended his work to the context of marked point processes and provided asymptotic normality by adding the assumption of finite range. Our goal is to deal with most classical models (see Baddeley and Turner (2000), Møller and Waggepetersen (2003) and Bertin et al. (1999b)) that could be interesting for practical purposes. They have been put into three categories according to their validity with respect to the previous existing works: 1) Overlap area point process 2) Multi-Strauss marked model, Strauss disc type process. 3) area-interaction, Geyer's triplet process, $k$-nearest-neighbour multi-Strauss marked model. Let us notice first that all the examples belong to the exponential family and satisfy the local stability and finite range properties. Due to the parametrization proposed by Jensen and Künsch (1994) and Mase (2000), they only consider the first category. By considering the exponential family, Jensen and Møller (1991) include a larger class of models. However, the required inhibition or hard-core type assumptions are only satisfied for examples 1 and 2. Examples 3 are only locally stable. What remains to be established is consistency for examples 3 and asymptotic normality for examples 2 and 3. In this paper, a general framework is proposed taking into consideration the previous remarks. Results are obtained using the general theory on minimum contrast estimators, *e.g.* Guyon (1995).

Section 2 introduces some background on marked Gibbs point processes. Our models are defined in Section 3. In the same spirit as Bertin et al. (1999a), providing existence results of stationary Gibbs states, assumptions on these models are expressed in terms of the local energy function. We also describe examples of interest of this work. Section 4 presents the pseudolikelihood method and our main results requiring two additional assumptions. The first one is an identifiability condition ensuring strong consistency. The second one is related to the definiteness of the asymptotic covariance matrix of the maximum pseudolikelihood estimator. These two assumptions allow us to derive a practical asymptotic normality result. They are verified for all the considered examples in Section 5. It is the belief of the authors that these assumptions are not restrictive since they should be true for every well-parametrized model. Proofs have been postponed until Section 6.

## 2. Background on marked Gibbs point processes

For the sake of simplicity, the framework of this paper is restricted to two-dimensional marked Gibbs point processes. All the results must remain valid in the general $d$-dimensional ($d \geq 1$) case. Define $\mathcal{B}^2$ the Borel $\sigma$-algebra on $\mathbb{R}^2$, $\mathcal{B}_b^2$ the set of bounded Borel susbsets of $\mathbb{R}^2$ and $\lambda^2$ the Lebesgue measure on $\mathbb{R}^2$. Denote also by $\mathbb{M}$, $\mathcal{M}$ and $\lambda^{\mathrm{m}}$ the mark space and its corresponding $\sigma$-algebra and probability measure. Let $\mathbb{S} := \mathbb{R}^2 \times \mathbb{M}$, $\mathcal{B} := \mathcal{B}^2 \otimes \mathcal{M}$ and $\mu := \lambda^2 \otimes \lambda^{\mathrm{m}}$ denote respectively the state space and its corresponding $\sigma$-algebra and measure.

For shortness, let us denote $x^m = (x, m)$ for any $x \in \mathbb{R}^2$ and any mark $m \in \mathbb{M}$ and $|\Lambda| := \lambda^2(\Lambda)$ for any $\Lambda \in \mathcal{B}^2$. In addition, $|I|$ designates the number of elements of some countable set $I$, $\Lambda^c$ is the complementary of some set $\Lambda$ in $\mathbb{R}^2$ and $||\cdot||$ is the $\ell^2$-norm. Let us define for all $i = (i_1, i_2) \in \mathbb{Z}^2$, $d > 0$ and $\rho \geq 0$ $\Delta_i(d) := \left\{z \in \mathbb{R}^2, d\left(i_j - \frac{1}{2}\right) \leq z_j \leq d\left(i_j + \frac{1}{2}\right), j = 1, 2\right\}$ and $\mathbb{B}(i, \rho) := \{k \in \mathbb{Z}^2 : |k - i| \leq \rho\}$ with $|i| := max(|i_1|, |i_2|)$.

Let $\widetilde{\Omega}$ denote the set of so-called configurations -of marked points- $\varphi := \{x_i^{m_i}\}_{i \in I}$ where $I$ is a subset of $\mathbb{N}$ and $((x_i, m_i))_{i \in I}$ is a sequence of elements of $\mathbb{S}$. In particular, any element $\varphi \in \widetilde{\Omega}$ has the following representation $\varphi = \sum_{i \in I} \delta_{x_i^{m_i}}$ as an integer-valued measure on $\mathbb{S}$ such that for every $F \in \mathcal{B}_b^2$, $\varphi(F) \in \mathbb{N}$, where $\delta_{x^m}$ is the Dirac measure at some element $x^m \in \mathbb{S}$. The subset of $\widetilde{\Omega}$ with elements $\varphi$ satisfying $|\varphi| := \varphi(\mathbb{S}) < +\infty$ is denoted by $\widetilde{\Omega}_f$. The space $\widetilde{\Omega}$ is equipped with the $\sigma$-algebra $\mathcal{F}$ generated by the family of sets $\left\{\varphi \in \widetilde{\Omega} : \varphi(F) = n\right\}$ with $n \in \mathbb{N}$ and $F \in \mathcal{B}_b^2$. For every $F \in \mathcal{B}^2$ and $\varphi \in \widetilde{\Omega}$ represented as $\varphi = \sum_{i \in I} \delta_{x_i^{m_i}}$, one introduces $\varphi_F := \sum_{i \in I, x_i^{m_i} \in F} \delta_{x_i^{m_i}}$ which can be viewed as the configuration of marked points of $\varphi$ in $F$. Furthermore, for every $\Lambda \in \mathcal{B}_b^2$, $\varphi_\Lambda$ conveniently denotes $\varphi_{\Lambda \times \mathbb{M}}$.

A marked point process is a $\widetilde{\Omega}$-valued random variable, denoted by $\Phi$, with probability distribution $P$ on $(\widetilde{\Omega}, \mathcal{F})$. The intensity measure $N_P$ of $P$ is defined as a measure on $\mathcal{B}^2$ such that for any $F \in \mathcal{B}_b^2$:

$$N_P(F) = \int_{\widetilde{\Omega}} \varphi(F) P(d\varphi) := \boldsymbol{E}(\Phi(F)).$$

In the stationary case, $N_P(F) = \nu_P \lambda^2(F)$ where $\nu_P$ is called the intensity of $P$. A marked Gibbs point process is usually defined using a family of local specifications with respect to a weight process (often a stationary marked Poisson process with distribution $Q$ and intensity $\lambda_Q = 1$). Let $\Lambda$ be a bounded region in $\mathbb{R}^2$. For such a process, given some configuration $\varphi_{\Lambda^c}$ on $\Lambda^c$, the conditional probability on $\Lambda$ is of the form, for any $F \in \mathcal{F}$:

$$\Pi_\Lambda(\varphi, F) = \left\{ \frac{1}{Z_\Lambda(\varphi)} \int_{\widetilde{\Omega}_\Lambda} e^{-V(\psi|\varphi_{\Lambda^c})} 1_F(\psi \cup \varphi_{\Lambda^c}) Q_\Lambda(d\psi) \right\} 1_{R_\Lambda}(\varphi),$$

with the partition function

$$Z_\Lambda(\varphi) = \int_{\widetilde{\Omega}_\Lambda} e^{-V(\psi|\varphi_{\Lambda^c})} Q_\Lambda(d\psi)$$

and $R_\Lambda = \{\varphi \in \widetilde{\Omega} : 0 < Z_\Lambda(\varphi) < +\infty\}$ where

$$\int f(\psi) Q_\Lambda(\psi) := e^{-\mu(\Lambda \times \mathbb{M})} \sum_{n=0}^{+\infty} \frac{1}{n!} \int f(\underbrace{\{x_1^{m_1}, \ldots, x_n^{m_n}\}}_{\psi}) d\mu^{\otimes n}(x_1^{m_1}, \ldots, x_n^{m_n}).$$

Let us define the subset of all admissible configurations

$$\Omega := \left\{ \varphi \in \widetilde{\Omega} : \varphi \in \cap_{\Lambda \in \mathcal{B}_b^2} R_\Lambda \right\}$$

and denote by $\Omega_f := \widetilde{\Omega}_f \cap \Omega$. Whereas the finite energy function $V(\varphi)$ (for any $\varphi \in \Omega_f$) measures the cost of any configuration, the local energy $V(\psi|\varphi)$ (for any $\varphi, \psi \in \Omega_f$) represents the energy required to add the points of $\psi$ in $\varphi$:

$$V(\psi|\varphi) = V(\psi \cup \varphi) - V(\varphi).$$

Let us notice that when $\psi$ is a singleton $\{x^m\}$, we denote by a slight abuse $V(x^m|\varphi)$ instead of $V(\{x^m\}|\varphi)$. It is well-known that the collection of probability kernels $(\Pi_\Lambda)_{\Lambda \in \mathcal{B}_b^2}$ satisfies the set of compatibility and measurability conditions which define a local specification in the Preston's sense (Preston (1976)). The main condition is the consistency:

$$\Pi_\Lambda \Pi_{\Lambda'} = \Pi_\Lambda \quad \text{for} \quad \Lambda' \subset \Lambda.$$

Notice that some conditions are needed to ensure the existence of a probability measure $P$ related to any local energy $V$ and any weight process that satisfies the so-called Dobrushin-Lanford-Ruelle (D.L.R.) equations:

$$P(F|\mathcal{F}_{\Lambda^c})(\varphi) = \Pi_\Lambda(\varphi, F) \quad \text{for } P \text{ a.e. } \varphi \in \Omega \quad \text{for any } \Lambda \in \mathcal{B}_b^2 \text{ and } F \in \mathcal{F}.$$

For the general theory of Gibbs point processes, the reader may refer to Kallenberg (1983); Stoyan et al. (1995) and the references therein.

For some finite configuration $\varphi$ (resp. some set $G$) and for all $x \in \mathbb{R}^2$, $\varphi_x$ (resp. $G_x$) denotes the configuration $\varphi$ (resp. the set $G$) translated of $x$. Finally, in this work a non-marked point process can be viewed as a particular case of marked point processes with $\mathbb{M} = \{0\}$.

## 3. Definitions and examples of marked Gibbs models

The framework of this paper is restricted to stationary marked Gibbs point processes based on an energy function invariant by translation, $V(\varphi;\boldsymbol{\theta})$, parametrized by some $\boldsymbol{\theta}\in\boldsymbol{\Theta}$, where $\boldsymbol{\Theta}$ is some compact set of $\mathbb{R}^p$. The model is also assumed to belong to an exponential family, i.e.

$$V(\varphi;\boldsymbol{\theta})=\boldsymbol{\theta}^T\boldsymbol{v}(\varphi), \tag{3.1}$$

where $\boldsymbol{v}(\varphi)=(v_1(\varphi),\ldots,v_p(\varphi))$ is the vector of sufficient statistics. The local energy is then expressed as

$$V(x^m|\varphi;\boldsymbol{\theta})=\boldsymbol{\theta}^T\boldsymbol{v}\left(x^m|\varphi\right), \tag{3.2}$$

where $\boldsymbol{v}(x^m|\varphi)=(v_1(x^m|\varphi),\ldots,v_p(x^m|\varphi)):=\boldsymbol{v}(\varphi\cup\{x^m\})-\boldsymbol{v}(\varphi)$.
Our models satisfy the general condition **[Mod]** described by the following statements:

**[Mod:S]** *Stability of the local energy*: there exists $K\geq 0$ such that for all $(m,\varphi)\in\mathbb{M}\times\Omega_f$

$$V(0^m|\varphi;\boldsymbol{\theta})\geq -K.$$

**[Mod:L]** *Locality of the local energy*: there exists $D\geq 0$ such that for all $(m,\varphi)\in\mathbb{M}\times\Omega_f$

$$V(0^m|\varphi;\boldsymbol{\theta})=V\left(0^m|\varphi_{B(0,D)};\boldsymbol{\theta}\right),$$

where $\mathcal{B}(x,r)$ denotes the ball centered at $x\in\mathbb{R}^2$ with radius $r>0$.

**[Mod:I]** *Integrability condition*: for $i=1,\ldots,p$, there exist $\kappa_i^{(\mathrm{sup})}\geq 0$, $k_i\in\mathbb{N}$ such that for all $(m,\varphi)\in\mathbb{M}\times\Omega_f$

$$v_i(0^m|\varphi)\leq\kappa_i^{(\mathrm{sup})}|\varphi_{\mathcal{B}(0,D)}|^{k_i}.$$

Let us notice that, unlike **[Mod:I]**, the assumptions **[Mod:S]** and **[Mod:L]** cannot be directly expressed only in terms of the sufficient statistics. Nevertheless, **[Mod]** is satisfied as soon as for $i=1,\ldots,p$, there exist $\kappa_i^{(\mathrm{inf})},\kappa_i^{(\mathrm{sup})}\geq 0$, $k_i\in\mathbb{N}$ such that one of both following assumptions is satisfied for all $(m,\varphi)\in\mathbb{M}\times\Omega_f$:

**[Mod-1]**

$$\theta_i\geq 0 \text{ and } -\kappa_i^{(\mathrm{inf})}\leq v_i(0^m|\varphi)=v_i(0^m|\varphi_{\mathcal{B}(0,D)})\leq\kappa_i^{(\mathrm{sup})}|\varphi_{\mathcal{B}(0,D)}|^{k_i}.$$

**[Mod-2]**

$$-\kappa_i^{(\mathrm{inf})}\leq v_i(0^m|\varphi)=v_i(0^m|\varphi_{\mathcal{B}(0,D)})\leq\kappa_i^{(\mathrm{sup})}.$$

Indeed, let $I_1$ and $I_2$ be the partition of $\{1,\ldots,p\}$ such that for any $i\in I_1$ (resp. $i\in I_2$), $v_i$ satisfies **[Mod-1]** (resp. **[Mod-2]**) then

$$\begin{aligned}
V(0^m|\varphi;\boldsymbol{\theta}) &= \sum_{i\in I_1}\theta_i v_i(0^m|\varphi)+\sum_{i\in I_2}\theta_i v_i(0^m|\varphi)\\
&\geq -p\left(\max_{i\in I_1}\Big(\theta_i\kappa_i^{(\mathrm{inf})}\Big)-\max_{i\in I_2}\Big(|\theta_i|\times\max(\kappa_i^{(\mathrm{inf})},\kappa_i^{(\mathrm{sup})})\Big)\right):=-K(\boldsymbol{\theta})
\end{aligned}$$

which ensures [**Mod:S**] with $K := \sup_{\boldsymbol{\theta}\in\boldsymbol{\Theta}} K(\boldsymbol{\theta})$ ([**Mod:L**] and [**Mod:I**] are clearly satisfied).
Let us also point out that

- the well-known characteristics [**Mod:S**] and [**Mod:L**] associated with finite energies that are translation invariant ensure the existence of stationary measures (see Bertin et al. (1999a)).
- the local stability implies the Ruelle-bound correlation function (see Ruelle (1970)) leading to: $\boldsymbol{E}\left(\left|\Phi_{\mathcal{B}(0,D)}\right|^{\beta}\right) < +\infty$ for any $\beta > 0$.
- A key-ingredient of our proofs is the following one: for any $\alpha > 0$, any $\boldsymbol{\theta} \in \boldsymbol{\Theta}$ and any $i = 1, \dots, p$

$$\boldsymbol{E}\left(|v_i(0^M|\Phi)|^{\alpha} e^{-\boldsymbol{\theta}^T \boldsymbol{v}\left(0^M|\Phi\right)}\right) < +\infty. \tag{3.3}$$

  This condition is fulfilled under the assumptions [**Mod:S**] and [**Mod:I**].

Let us now present some examples. Except the one based on the $k$-nearest-neighbour graph, all examples are really classical and can be found *e.g.* in Baddeley and Turner (2000) and Møller and Waggepetersen (2003). For each example, we present the model through the sufficient statistics, the set of the parameters values (including $\boldsymbol{\Theta}$) for which the model is defined in the litterature and then we verify that [**Mod**] is satisfied. This proves in particular the existence of stationary Gibbs states in $\mathbb{R}^2$.

First of all, note that when $v_i(0^m|\varphi) := 1$, then $v_i$ obviously satisfies [**Mod-2**]. Recall that for a non-marked point process $\mathbb{M} = \{0\}$.

*Overlap area point process*

- This non-marked process is defined for $p = 2$ and $R > 0$ by

$$v_1(\varphi) := |\varphi| \text{ and } v_2(\varphi) := \sum_{\{x,y\}\in\mathcal{P}_2(\varphi)} \left|\mathcal{B}(x, R/2) \cap \mathcal{B}(y, R/2)\right|,$$

  where $\mathcal{P}_k(\varphi)$ ($k \geq 1$) is the set of all subsets of $\varphi$ with $k$ elements. Alternatively,

$$v_1(0|\varphi) := 1 \text{ and } v_2(0|\varphi) := \sum_{x\in\varphi} \left|\mathcal{B}(0, R/2) \cap \mathcal{B}(x, R/2)\right|.$$

  Let us notice that

$$\begin{aligned}&\left|\mathcal{B}(0, R/2) \cap \mathcal{B}(x, R/2)\right| \\ &\quad = \frac{1}{2}\left(R^2 Arcos\frac{||x||}{R} - ||x||\sqrt{R^2 - ||x||^2}\right) \mathbf{1}_{[0,R]}(||x||).\end{aligned}$$

- $\boldsymbol{\theta} \in \mathbb{R} \times \mathbb{R}^+$.
- $v_2$ satisfies [**Mod-1**] with $\kappa_2^{(\text{inf})} = 0$, $\kappa_2^{(\text{sup})} = \frac{\pi R^2}{4}$, $D = R$ and $k_2 = 1$.

*Multi-Strauss marked point process*

- Let $\mathbb{M} = \{1, \ldots, M\}$, $\lambda^{\mathrm{m}}$ the uniform probability measure on $\mathbb{M}$ and $p \geq 2$. Decompose $\varphi = \cup_{m=1}^{M} \varphi^m$ with $\varphi^m := \{x^m \in \Omega : x^m \in \varphi\}$ for any $m \in \mathbb{M}$. The finite energy of this process is defined by

$$\begin{aligned} V(\varphi;\boldsymbol{\theta}) \quad &:= \quad \sum_{m_1=1}^{M} \theta_1^{m_1,m_1} v_1^{m_1,m_1}(\varphi^{m_1}) \\ &+ \sum_{1\leq m_1\leq m_2\leq M} \sum_{i=2}^{p^{m_1,m_2}} \theta_i^{m_1,m_2} v_i^{m_1,m_2}(\varphi^{m_1}\cup\varphi^{m_2}) \end{aligned}$$

with for all $m_1, m_2 \in \mathbb{M}$ and $i = 2, \ldots, p$

$$\begin{aligned} v_1^{m_1,m_1}(\varphi) \quad &:= \quad v_1^{m_1,m_1}(\varphi^{m_1}) := |\varphi^{m_1}| \\ v_i^{m_1,m_2}(\varphi) \quad &:= \quad v_i^{m_1,m_2}(\varphi^{m_1}\cup\varphi^{m_2}) \\ &= \quad \sum_{\{x_1^{m_1},x_2^{m_2}\}\in\mathcal{P}_2(\varphi)} \mathbf{1}_{[D_{i-1}^{m_1,m_2},D_i^{m_1,m_2}[}(||x_1-x_2||) \end{aligned} \tag{3.4}$$

where $0 \leq D_1^{m_1,m_2} < D_2^{m_1,m_2} < \ldots < D_{p^{m_1,m_2}}^{m_1,m_2} < +\infty$. In particular, the vector $\boldsymbol{\theta}$ could be ordered as follows:

$$\boldsymbol{\theta} = (\boldsymbol{\theta}^1, \boldsymbol{\theta}^2, \ldots, \boldsymbol{\theta}^M) \text{ where } \boldsymbol{\theta}^m = (\boldsymbol{\theta}^{m,m}, \boldsymbol{\theta}^{m,m+1}, \ldots, \boldsymbol{\theta}^{m,M})$$

with

$$\boldsymbol{\theta}^{m_1,m_2} = \begin{cases} (\theta_1^{m_1,m_2}, \theta_2^{m_1,m_2}, \ldots, \theta_{p^{m_1,m_2}}^{m_1,m_2}) & \text{if } m_1 = m_2 \\ (\theta_2^{m_1,m_2}, \theta_3^{m_1,m_2}, \ldots, \theta_{p^{m_1,m_2}}^{m_1,m_2}) & \text{otherwise.} \end{cases}$$

where $p = M + \sum_{1\leq m_1\leq m_2\leq M} (p^{m_1,m_2} - 1)$. One then derives the expression of the local energy

$$V(0^{m_1}|\varphi;\boldsymbol{\theta}) = \theta_1^{m_1,m_1} + \sum_{m_2=1}^{M} \sum_{i=2}^{p^{m_1,m_2}} \theta_i^{m_1,m_2} v_i^{m_1,m_2}(0^{m_1}|\varphi^{m_2})$$

where for convenience $\theta_i^{m_1,m_2}$ and $\theta_i^{m_2,m_1}$ denote the same parameter.
- For all $m_1, m_2 \in \mathbb{M}$ and $i = 2, \ldots, p^{m_1,m_2}$: $\theta_1^{m_1,m_1} \in \mathbb{R}$ and $\theta_i^{m_1,m_2} \in \begin{cases} \mathbb{R} & \text{when } D_1^{m_1,m_2} = \delta > 0 \\ \mathbb{R}^+ & \text{when } D_1^{m_1,m_2} = 0. \end{cases}$
- When $D_1^{m_1,m_2} = 0$ and $\theta_i^{m_1,m_2} \geq 0$, $v_i^{m_1,m_2}$ satisfies **[Mod-1]** with $\kappa_{i(m_1,m_2)}^{(\mathrm{inf})} = 0$, $\kappa_{i(m_1,m_2)}^{(\mathrm{sup})} = 1$, $D = \max_{m_1\leq m_2} D_{p^{m_1,m_2}}^{m_1,m_2}$ and $k_{i(m_1,m_2)} = 1$, where $\theta_i^{m_1,m_2} = \theta_{i(m_1,m_2)}$ is the $i(m_1,m_2)-$th element of the vector $\boldsymbol{\theta}$. Under the hard-core assumption $D_1^{m_1,m_2} = \delta > 0$, $v_i^{m_1,m_2}$ satisfies **[Mod-2]** with $\kappa_{i(m_1,m_2)}^{(\mathrm{inf})} = 0$, $\kappa_{i(m_1,m_2)}^{(\mathrm{sup})} = \lceil \frac{D^2}{\delta^2} \rceil$ and $D = \max_{m_1\leq m_2} D_{p^{m_1,m_2}}^{m_1,m_2}$.

*k−nearest-neighbour multi-Strauss marked point process*

- This marked point process is defined similarly as the multi-Strauss marked point process except that the complete graph $\mathcal{P}_2(\varphi)$ in (3.4) is replaced by the $k$-nearest-neighbour graph ($k \geq 1$).
- For all $m_1, m_2 \in \mathbb{M}$ and $i = 2, \ldots, p^{m_1,m_2}$: $\theta_i^{m_1,m_2} \in \mathbb{R}$.
- In Bertin et al. (1999b), it is proved that $v_i^{m_1,m_2}$ satisfies [**Mod-2**] with $\kappa_{i(m_1,m_2)}^{(\text{inf})} = \kappa_{i(m_1,m_2)}^{(\text{sup})} = 13k$ and $D = 2\max_{m_1 \leq m_2} D_{p^{m_1,m_2}}^{m_1,m_2}$.

*Strauss type disc process*

- Let $\mathbb{M} = [0, M_{\max}]$ with $0 < M_{\max} < +\infty$, $\lambda^{\mathrm{m}}$ the uniform probability measure on $\mathbb{M}$ and $p = 2$. This model is defined by

$$v_1(\varphi) = |\varphi| \text{ and } v_2(\varphi) = \sum_{\{x_1^{m_1}, x_2^{m_2}\} \in \mathcal{P}_2(\varphi)} \mathbf{1}_{[0, m_1+m_2]}(||x_2 - x_1||).$$

  Alternatively,

$$v_2(0^m|\varphi) = \sum_{x^{m'} \in \varphi} \mathbf{1}_{[0, m+m']}(||x||).$$

- $\boldsymbol{\theta} \in \mathbb{R} \times \mathbb{R}^+$.
- $v_2$ satisfies [**Mod-1**] with $\kappa_2^{(\text{inf})} = 0$, $\kappa_2^{(\text{sup})} = 1$, $D = M_{\max}$ and $k_2 = 1$.

*Geyer's triplet interaction point process*

- This non-marked point process is defined for $p = 3$ and $R > 0$ by

$$V(\varphi; \boldsymbol{\theta}) = \theta_1 |\varphi| + \theta_2 v_2(\varphi) + \theta_3 v_3(\varphi)$$

  where

$$v_2(\varphi) = \sum_{\{x_1, x_2\} \in \mathcal{P}_2(\varphi)} \mathbf{1}_{[0,R]}(||x_1 - x_2||)$$

  and

$$v_3(\varphi) = \sum_{\xi \in \mathcal{P}_3(\varphi)} \prod_{\{x_1, x_2\} \in \mathcal{P}_2(\xi)} \mathbf{1}_{[0,R]}(||x_1 - x_2||).$$

  Note that $v_3(\varphi)$ represents the number of triangles of $\varphi$ with edges of lengths lower than $R$. Alternatively

$$v_2(0|\varphi) = \sum_{x \in \varphi} \mathbf{1}_{[0,R]}(||x||)$$

  and

$$v_3(0|\varphi) = \sum_{\{x,y\} \in \mathcal{P}_2(\varphi)} \prod_{\{x_1, x_2\} \in \mathcal{P}_2(\{x,y,0\})} \mathbf{1}_{[0,R]}(||x_1 - x_2||).$$

- $\boldsymbol{\theta} \in \mathbb{R}^2 \times \mathbb{R}^+ \setminus \{0\}$.
- When $\theta_2 \geq 0$, $v_2$ and $v_3$ satisfy [**Mod-1**] with $D = R$, $\kappa_2^{(\text{inf})} = \kappa_3^{(\text{inf})} = 0$, $\kappa_2^{(\text{sup})} = \kappa_3^{(\text{sup})} = 1$, $k_2 = 1$ and $k_3 = 2$. When $\theta_2 < 0$, $v_2$ satisfies neither [**Mod-1**] nor [**Mod-2**]. However, Geyer (1999) proved that the local energy is stable and local (i.e. [**Mod:S**] and [**Mod:L**]) and [**Mod:I**] is satisfied with $D = R$, $\kappa_2^{(\text{sup})} = \kappa_3^{(\text{sup})} = 1$, $k_2 = 1$ and $k_3 = 2$.

*Area interaction point process*

- This model is the one-type marginal of the two-type Widom-Rowlinson model. Let $p = 2$ and $R > 0$

$$V(\varphi; \boldsymbol{\theta}) = \theta_1 |\varphi| + \theta_2 v_2(\varphi), \quad \text{with} \quad v_2(\varphi) := |\cup_{x \in \varphi} \mathcal{B}(x, R)|.$$

  Note that $v_2(\varphi)$ represents the area of the union of discs of radius $R$ centered at the points. Alternatively,

$$v_2(0|\varphi) := \left| \cup_{x \in \varphi_{\mathcal{B}(0,2R)} \cup \{0\}} \mathcal{B}(x, R) \setminus \cup_{x \in \varphi_{\mathcal{B}(0,2R)}} \mathcal{B}(x, R) \right|.$$

- $\boldsymbol{\theta} \in \mathbb{R}^2$.
- $v_2$ satisfies [**Mod-2**] with $\kappa_2^{(\text{inf})} = 0$, $\kappa_2^{(\text{sup})} = \pi D^2$ and $D = 2R$.

**Remark 1.** *Jensen and Møller (1991) have already proved the consistency property in the inhibition case, i.e. [**Mod:S**] with $K = 0$. In particular, they did not consider the area-interaction model with negative parameters $\theta_1$ and $\theta_2$ (for instance). This gap has now been filled by extending the result in the case when [**Mod:S**] is satisfied. However, unlike these authors, we require the additional assumption [**Mod:I**]. In order to simplify our assumptions, we deliberately decided not to propose this particular case since this last assumption is not restrictive and is satisfied for all examples considered above.*

**Remark 2.** *Note that unlike the Multi-Strauss marked point process, neither inhibition nor hard-core assumption is required for the k−nearest-neighbour multi-Strauss marked point process since its local energy is naturally stable.*

**Remark 3.** *Concerning the Geyer's triplet process, the case $\theta_3 = 0$ is not considered since it is a particular case of a multi-Strauss marked point process.*

**Remark 4.** *Through these different examples, one may note that some parameters are assumed to be known: for example, the parameters $D_i^{m_1,m_2}$ for the multi-Strauss marked point process, the hard-core parameter $\delta$, the parameter $M_{max}$ for the Strauss type disc process, the parameter $R$ for the Geyer's triplet process or the area interaction point process, .... Their estimations could be investigated by using ad hoc methods.*

## 4. MPLE: presentation and asymptotic results

### 4.1. Pseudolikelihood

As specified in the introduction, the idea of maximum pseudolikelihood is due to Besag (1974) who first introduced the concept for Markov random fields in order to avoid the normalizing constant. This work was then widely extended and Jensen and Møller (1991) (Theorem 2.2) obtained a general expression for marked Gibbs point processes. With our notation and up to a scalar factor, the pseudolikelihood defined for a configuration $\varphi$ and a domain of observation $\Lambda$ is denoted by $PL_\Lambda(\varphi;\boldsymbol{\theta})$ and given by

$$PL_\Lambda(\varphi;\boldsymbol{\theta}) = \exp\left(-\int_{\Lambda\times\mathbb{M}} e^{-V(x^m|\varphi;\boldsymbol{\theta})}\mu(dx^m)\right)\prod_{x^m\in\varphi_\Lambda} e^{-V(x^m|\varphi\setminus x^m;\boldsymbol{\theta})}. \quad (4.1)$$

It is more convenient to define (and work with) the log-pseudolikelihood, denoted by $LPL_\Lambda(\varphi;\boldsymbol{\theta})$.

$$LPL_\Lambda(\varphi;\boldsymbol{\theta}) = -\int_{\Lambda\times\mathbb{M}} e^{-V(x^m|\varphi;\boldsymbol{\theta})}\mu(dx^m) - \sum_{x^m\in\varphi_\Lambda} V(x^m|\varphi\setminus x^m;\boldsymbol{\theta}). \quad (4.2)$$

Our data consist in the realization of a point process with energy function $V(\cdot;\boldsymbol{\theta}^\star)$ satisfying **[Mod]**. Thus, $\boldsymbol{\theta}^\star$ is the true parameter to be estimated and it is assumed that $\boldsymbol{\theta}^\star \in \mathring{\boldsymbol{\Theta}}$. The Gibbs measure will be denoted by $P_{\boldsymbol{\theta}^\star}$. Moreover the point process is assumed to be observed in a domain $\Lambda_n \oplus D^\vee = \cup_{x\in\Lambda_n}\mathcal{B}(x,D^\vee)$ for some $D^\vee \geq D$. For the asymptotic normality result, it is also assumed that $\Lambda_n \subset \mathbb{R}^2$ can de decomposed into $\cup_{i\in I_n}\Delta_i$ where $I_n = \mathbb{B}(0,n)$ and for $i=(i_1,i_2)\in\mathbb{Z}^2$, $\Delta_i = \Delta_i(\widetilde{D})$ for some $\widetilde{D}>0$ fixed from now on. As a consequence, as $n\to+\infty$, $\Lambda_n\to\mathbb{R}^2$ such that $|\Lambda_n|\to+\infty$ and $\dfrac{|\partial\Lambda_n|}{|\Lambda_n|}\to 0$.

Define for any configuration $\varphi$, $U_n(\varphi;\boldsymbol{\theta}) = -\frac{1}{|\Lambda_n|}LPL_{\Lambda_n}(\varphi;\boldsymbol{\theta})$. The maximum pseudolikelihood estimate (MPLE) denoted by $\widehat{\boldsymbol{\theta}}_n(\varphi)$ is then defined by

$$\widehat{\boldsymbol{\theta}}_n(\varphi) = \operatorname*{arg\,max}_{\boldsymbol{\theta}\in\boldsymbol{\Theta}} LPL_{\Lambda_n}(\varphi;\boldsymbol{\theta}) = \operatorname*{arg\,min}_{\boldsymbol{\theta}\in\boldsymbol{\Theta}} U_n(\varphi;\boldsymbol{\theta}).$$

We will also need the following basic notations:

- Gradient vector of $U_n$: $\boldsymbol{U}_n^{(1)}(\varphi;\boldsymbol{\theta}) := -|\Lambda_n|^{-1}\boldsymbol{LPL}_{\Lambda_n}^{(1)}(\varphi;\boldsymbol{\theta})$ where for any bounded Borel set $\Lambda$, $\left(\boldsymbol{LPL}_\Lambda^{(1)}(\varphi;\boldsymbol{\theta})\right)_j$ is defined for $j=1,\ldots,p$ by

$$\left(\boldsymbol{LPL}_\Lambda^{(1)}(\varphi;\boldsymbol{\theta})\right)_j = \int_{\Lambda\times\mathbb{M}} v_j(x^m|\varphi)e^{-V(x^m|\varphi;\boldsymbol{\theta})}\mu(dx^m) - \sum_{x^m\in\varphi_\Lambda} v_j(x^m|\varphi\setminus x^m)$$

- Hessian matrix of $U_n$: $\underline{\boldsymbol{U}}_n^{(2)}(\varphi;\boldsymbol{\theta}) := -|\Lambda_n|^{-1}\underline{\boldsymbol{LPL}}_{\Lambda_n}^{(2)}(\varphi;\boldsymbol{\theta})$ where for any bounded Borel set $\Lambda$, $\left(\underline{\boldsymbol{LPL}}_{\Lambda}^{(2)}(\varphi;\boldsymbol{\theta})\right)_{j,k}$ is defined for $j,k=1,\ldots,p$ by

$$\left(\underline{\boldsymbol{LPL}}_{\Lambda}^{(2)}(\varphi;\boldsymbol{\theta})\right)_{j,k} = \int_{\Lambda\times\mathbb{M}} v_j(x^m|\varphi)v_k(x^m|\varphi)e^{-V(x^m|\varphi;\boldsymbol{\theta})}\mu(dx^m)$$

Finally, note that from the decomposition of the observation domain $\Lambda_n$, one has

$$\boldsymbol{U}_n^{(1)}(\varphi;\boldsymbol{\theta}) = |\Lambda_n|^{-1}\sum_{i\in I_n}\boldsymbol{LPL}_{\Delta_i}^{(1)}(\varphi;\boldsymbol{\theta})$$

and

$$\underline{\boldsymbol{U}}_n^{(2)}(\varphi;\boldsymbol{\theta}) = |\Lambda_n|^{-1}\sum_{i\in I_n}\underline{\boldsymbol{LPL}}_{\Delta_i}^{(2)}(\varphi;\boldsymbol{\theta}).$$

### *4.2. Asymptotic results of the MPLE*

This section provides consistency and asymptotic normality of the maximum pseudolikelihood estimator. Let us first consider the following assumption

**[Ident]** Identifiability condition: there exists $A_1,\ldots,A_\ell$, $\ell\geq p$ events of $\Omega$ and $A_1^{\mathbb{m}},\ldots,A_\ell^{\mathbb{m}}$ events of $\mathcal{M}$ such that:

- the $\ell$ events $B_j := A_j^{\mathbb{m}}\times A_j$ are disjoint and satisfy $\lambda^{\mathbb{m}}\otimes P_{\boldsymbol{\theta}^\star}(B_j)>0$
- for all $((m_1,\varphi_1),\ldots,(m_\ell,\varphi_\ell))\in B_1\times\cdots\times B_\ell$ the $(\ell,p)$ matrix with entries $v_j(0^{m_i}|\varphi_i)$ is injective.

**Theorem 1.** *Under the assumptions* **[Mod]** *and* **[Ident]***, for $P_{\boldsymbol{\theta}^\star}-$almost every $\varphi$, the maximum pseudolikelihood estimate $\widehat{\boldsymbol{\theta}}_n(\varphi)$ converges towards $\boldsymbol{\theta}^\star$ as $n$ tends to infinity.*

For the next result consider

**[SDP]** For some $\overline{\Lambda} := \cup_{i\in\mathbb{B}\left(0,\left\lceil\frac{D}{\overline{D}}\right\rceil\right)}\Delta_i(\overline{D})$ with $\overline{D}>0$, there exists $A_0,\ldots,A_\ell$, $\ell\geq p$ disjoint events of $\overline{\Omega} := \left\{\varphi\in\Omega : \varphi_{\Delta_i(\overline{D})}=\emptyset, 1\leq|i|\leq 2\left\lceil\frac{D}{\overline{D}}\right\rceil\right\}$ such that

- for $j=0,\ldots,\ell$, $P_{\boldsymbol{\theta}^\star}(A_j)>0$.
- for all $(\varphi_0,\ldots,\varphi_\ell)\in A_0\times\cdots\times A_\ell$ the $(\ell,p)$ matrix with entries $\left(\boldsymbol{LPL}_{\overline{\Lambda}}^{(1)}(\varphi_i;\boldsymbol{\theta})\right)_j - \left(\boldsymbol{LPL}_{\overline{\Lambda}}^{(1)}(\varphi_0;\boldsymbol{\theta})\right)_j$ is injective.

**Theorem 2.** *Under the assumptions* **[Mod]** *and* **[Ident]***, we have, for any fixed $\widetilde{D}$, the following convergence in distribution as $n\to+\infty$*

$$|\Lambda_n|^{1/2}\,\underline{\boldsymbol{U}}_n^{(2)}(\Phi;\boldsymbol{\theta}^\star)\,\left(\widehat{\boldsymbol{\theta}}_n(\Phi)-\boldsymbol{\theta}^\star\right)\to\mathcal{N}\left(0,\underline{\boldsymbol{\Sigma}}(\boldsymbol{\theta}^\star)\right), \tag{4.3}$$

*where*

$$\underline{\boldsymbol{\Sigma}}(\boldsymbol{\theta}^\star) = \sum_{i \in \mathbb{B}(0,\lceil D \rceil)} \boldsymbol{E}\left( \boldsymbol{LPL}^{(1)}_{\Delta_0(1)}(\Phi;\boldsymbol{\theta}^\star)\, \boldsymbol{LPL}^{(1)}_{\Delta_i(1)}(\Phi;\boldsymbol{\theta}^\star)^T \right). \tag{4.4}$$

*In addition under the assumption* ***[SDP]***

$$|\Lambda_n|^{1/2}\, \widehat{\underline{\boldsymbol{\Sigma}}}_n(\Phi; D^\vee, \widetilde{D}, \widehat{\boldsymbol{\theta}}_n(\Phi))^{-1/2}\, \underline{\boldsymbol{U}}_n^{(2)}(\Phi; \widehat{\boldsymbol{\theta}}_n(\Phi)) \left( \widehat{\boldsymbol{\theta}}_n(\Phi) - \boldsymbol{\theta}^\star \right) \to \mathcal{N}\left(0, \underline{\boldsymbol{I}}_p\right), \tag{4.5}$$

*where for some* $\boldsymbol{\theta}$ *and any configuration* $\varphi$*, the matrix* $\widehat{\underline{\boldsymbol{\Sigma}}}_n(\varphi; D^\vee, \widetilde{D}, \boldsymbol{\theta})$ *is defined by*

$$\widehat{\underline{\boldsymbol{\Sigma}}}_n(\varphi; D^\vee, \widetilde{D}, \boldsymbol{\theta}) = |\Lambda_n|^{-1} \sum_{i \in I_n} \sum_{j \in \mathbb{B}\left(i, \left\lceil \frac{D^\vee}{\widetilde{D}} \right\rceil\right) \cap I_n} \boldsymbol{LPL}^{(1)}_{\Delta_i}(\varphi;\boldsymbol{\theta})\, \boldsymbol{LPL}^{(1)}_{\Delta_j}(\varphi;\boldsymbol{\theta})^T. \tag{4.6}$$

**Remark 5.** *Let us underline that the scaling that yields to asymptotic normality corresponds to the usual parametric rate. Indeed, in the d-dimensional case one would obtain*

$$|\Lambda_n|^{1/2} = |\Lambda_n(\widetilde{D})|^{1/2} = |I_n|^{1/2} \times |\Delta_i(\widetilde{D})|^{1/2} = \widetilde{D}^{d/2}((2n+1)^d)^{1/2} \sim (2\widetilde{D})^{d/2}(n^d)^{1/2}.$$

**Remark 6.** *We would like to underline that* $\underline{\boldsymbol{\Sigma}}(\boldsymbol{\theta}^\star) = \underline{\boldsymbol{\Sigma}}(\widetilde{D}, \boldsymbol{\theta}^\star) = \underline{\boldsymbol{\Sigma}}(1, \boldsymbol{\theta}^\star)$ *where for all* $\widetilde{D} > 0$

$$\underline{\boldsymbol{\Sigma}}(\widetilde{D}, \boldsymbol{\theta}^\star) = \widetilde{D}^{-2} \sum_{i \in \mathbb{B}\left(0, \left\lceil \frac{D}{\widetilde{D}} \right\rceil\right)} \boldsymbol{E}\left( \boldsymbol{LPL}^{(1)}_{\Delta_0}(\Phi;\boldsymbol{\theta}^\star)\, \boldsymbol{LPL}^{(1)}_{\Delta_i}(\Phi;\boldsymbol{\theta}^\star)^T \right). \tag{4.7}$$

## 5. Back to examples

This section is devoted to proving that all of our examples satisfy both assumptions **[Ident]** and **[SDP]**. For the assumption **[Ident]** $\underline{\boldsymbol{V}}$ denotes the matrix with entries $v_j(0^{m_i}|\varphi_i)$ where $(m_i, \varphi_i) \in B_i$ have to be defined according to the different examples. The assumption **[SDP]** may be rewritten for all $k = 1, \ldots, \ell$ and for all $\varphi_k \in A_k$ and $\varphi_0 \in A_0$:

$$\left( \forall \boldsymbol{y} \in \mathbb{R}^p, \boldsymbol{y}^T \left( \boldsymbol{LPL}^{(1)}_{\overline{\Lambda}}(\varphi_k;\boldsymbol{\theta}^\star) - \boldsymbol{LPL}^{(1)}_{\overline{\Lambda}}(\varphi_0;\boldsymbol{\theta}^\star) \right) = \boldsymbol{y}^T (\boldsymbol{L}(\varphi_k;\boldsymbol{\theta}^\star) - \boldsymbol{R}(\varphi_k)) = 0 \right) \Rightarrow \boldsymbol{y} = 0,$$

where for any configuration $\varphi \in \overline{\Omega}$ and $\varphi_0 \in A_0$

$$\begin{aligned} \boldsymbol{L}(\varphi;\boldsymbol{\theta}^\star) &:= \int_{\overline{\Lambda} \times \mathbb{M}} \boldsymbol{v}(x^m|\varphi)\, e^{-\boldsymbol{\theta}^{\star T} \boldsymbol{v}(x^m|\varphi)} \mu(dx^m) \\ &\quad - \int_{\overline{\Lambda} \times \mathbb{M}} \boldsymbol{v}(x^m|\varphi_0)\, e^{-\boldsymbol{\theta}^{\star T} \boldsymbol{v}(x^m|\varphi_0)} \mu(dx^m) \\ \boldsymbol{R}(\varphi) &:= \sum_{x^m \in \varphi \cap \overline{\Lambda}} \boldsymbol{v}(x^m|\varphi \setminus x^m) - \sum_{x^m \in \varphi_0 \cap \overline{\Lambda}} \boldsymbol{v}(x^m|\varphi_0 \setminus x^m). \end{aligned}$$

Concerning this assumption, we choose $\overline{D} > D$ in all our examples.

### 5.1. Overlap area point process

*Assumption **[Ident]***

Consider

$$\begin{aligned} A_1 &:= \left\{\varphi \in \Omega : \varphi_{\mathcal{B}(0,D)} = \emptyset\right\} \\ A_2 &:= \left\{\varphi \in \Omega : \varphi_{\mathcal{B}(0,D)} = \{z\}, z \in \mathcal{B}((0,D/2),D/4)\right\}. \end{aligned}$$

We have for all $(\varphi_1, \varphi_2) \in A_1 \times A_2$

$$\underline{\boldsymbol{V}} = \begin{pmatrix} 1 & 0 \\ 1 & v_2(0|\varphi_2) \end{pmatrix}.$$

For every $\varphi_2 \in A_2$ such that $(\varphi_2)_{\mathcal{B}(0,D)} = \{z\}$ with for all $z \in \mathcal{B}((0,D/2),D/4)$, one remarks that $|\mathcal{B}(0,R/2) \cap \mathcal{B}(z,R/2)| := g_2(||z||) > 0$, then $\det(\underline{\boldsymbol{V}}) \neq 0$.

*Assumption **[SDP]***

Denote by $A_0$ any configuration set. Consider $A_n(\eta)$ for $n \geq 1$ and for some $0 < \eta < D$ the following configuration set

$$A_n(\eta) = \left\{\varphi \in \overline{\Omega} : \varphi_{\Delta_0(\overline{D})} = \{z_1, \ldots, z_n\} \text{ with } z_1, \ldots, z_n \in \mathcal{B}(0,\eta)\right\}.$$

For any $\varphi_n \in A_n(\eta)$, we have

$$\left| \int_{\overline{\Lambda}} v_j(x|\varphi_n) e^{-\boldsymbol{\theta}^{\star T} \boldsymbol{v}(x|\varphi_n)} dx \right| \leq \begin{cases} |\overline{\Lambda}| e^K & \text{if j=1} \\ n\pi R^2 |\overline{\Lambda}| e^K & \text{if j=2.} \end{cases}$$

where $K$ comes from the local stability property. Let us also remark that

$$\boldsymbol{y}^T \boldsymbol{R}(\varphi_n) = n y_1 + y_2 \times \sum_{x \in \varphi_n \cap \overline{\Lambda}} v_2(x|\varphi_n \setminus x) - \boldsymbol{y}^T \sum_{x \in \varphi_0 \cap \overline{\Lambda}} \boldsymbol{v}(x|\varphi_0 \setminus x)$$

with

$$0 < n(n-1) g_2(\eta) \leq \sum_{x \in \varphi_n \cap \overline{\Lambda}} v_2(x|\varphi_n \setminus x) \leq n(n-1) \frac{\pi R^2}{4}.$$

Therefore by combining these arguments, for every $\varepsilon > 0$, we have for $n$ large enough

$$\begin{aligned} |y_2| g_2(\eta) &\leq \left| \frac{1}{n(n-1)} y_2 \sum_{x \in \varphi_n \cap \overline{\Lambda}} v_2(x|\varphi_n \setminus x) \right| \\ &= \left| \frac{1}{n(n-1)} \left( y_2 \sum_{x \in \varphi_n \cap \overline{\Lambda}} v_2(x|\varphi_n \setminus x) + \boldsymbol{y}^T (\boldsymbol{L}(\varphi_n; \boldsymbol{\theta}^\star) - \boldsymbol{R}(\varphi_n)) \right) \right| \leq \varepsilon. \end{aligned}$$

By choosing $\varepsilon = \frac{|y_2| g_2(\eta)}{2}$, this leads to $y_2 = 0$. Then, for every $\varepsilon' > 0$ we may obtain for $n$ large enough

$$|y_1| = \left| y_1 + \frac{1}{n}(y_1, 0)^T \left( \boldsymbol{L}(\varphi_n; \boldsymbol{\theta}^\star) - \boldsymbol{R}(\varphi_n) \right) \right| \le \varepsilon'.$$

By choosing $\varepsilon' = |y_1|/2$, this leads to $y_1 = 0$.

### *5.2. Multi-Strauss marked type models*

*Assumption* ***[Ident]***

Define for any $m, m_1, m_2 \in \{1, \ldots, M\}$ with $m_2 \ge m_1$ and for any $i = 2, \ldots, p^{m_1,m_2}$

$$\begin{aligned} A_0 &:= \left\{ \varphi \in \Omega : \varphi_{\mathcal{B}(0,D)} = \emptyset \right\} \\ A_m^{\mathbb{m}} &:= \{m\} \\ A_i^{m_1,m_2} &:= \left\{ \varphi \in \Omega : \varphi_{\mathcal{B}(0,D)} = \{z^{m_1}\}, \text{ with } z \in \mathcal{B}(0, D_i^{m_1,m_2}) \setminus \mathcal{B}(0, D_{i-1}^{m_1,m_2}) \right\} \end{aligned}$$

The following events $B_i^{m_1,m_2}$ are defined for any $i = 1, \ldots, p^{m_1,m_1}$ when $m_1 = m_2$ and any $i = 2, \ldots, p^{m_1,m_2}$ when $m_1 < m_2$ such that:

$$B_i^{m_1,m_2} = \begin{cases} A_{m_1}^{\mathbb{m}} \times A_0 & \text{if } m_1 = m_2 \text{ and } i = 1 \\ A_{m_1}^{\mathbb{m}} \times A_i^{m_1,m_2} & \text{otherwise} \end{cases}$$

One may order these $\ell = p$ events as $B_1, \ldots, B_\ell$ where $B_{k_i^{m_1,m_2}} := B_i^{m_1,m_2}$ with

$$k_i^{m_1,m_2} = (m_1 - 1 + \delta_{m_1,m_2}) + \left( i - 1 + \sum_{m_2'=m_1}^{m_2-1} (p^{m_1,m_2'} - 1) + \sum_{m_1'=1}^{m_1-1} \sum_{m_2=m_1'}^{M} (p^{m_1',m_2} - 1) \right).$$

The corresponding matrix is then $\underline{\boldsymbol{V}} = \begin{pmatrix} V_1 & 0 & \cdots & 0 \\ 0 & \ddots & \ddots & \vdots \\ \vdots & \ddots & \ddots & 0 \\ 0 & \cdots & 0 & V_M \end{pmatrix}$ with

$$V_m = \begin{pmatrix} 1 & V_{m,m} & 0 & \cdots & 0 \\ \vdots & 0 & \ddots & & \vdots \\ \vdots & \vdots & \ddots & \ddots & 0 \\ 1 & 0 & \cdots & 0 & V_{m,M} \end{pmatrix} \text{ and } V_{m_1,m_2} = \begin{pmatrix} 0 & \cdots & \cdots & 0 \\ 1 & \ddots & & \vdots \\ 0 & \ddots & \ddots & \vdots \\ \vdots & \ddots & \ddots & 0 \\ 0 & \cdots & 0 & 1 \end{pmatrix}.$$

*Assumption **[SDP]***

Let us first introduce the following sets for any $\eta > 0$ and $d > 0$

$$A(\eta, d) = \left\{ (z_1, z_2) \in \Delta_0(\overline{D})^2 : z_1 \in \mathcal{B}\left((0,0), \frac{\eta}{4}\right) \text{ and } z_2 \in \mathcal{B}\left((d,0), \frac{3\eta}{4}\right)\right\}$$

$$A_-(\eta, d) = \left\{ (z_1, z_2) \in \Delta_0(\overline{D})^2 : z_1 \in \mathcal{B}\left((0,0), \frac{\eta}{4}\right) \text{ and } z_2 \in \mathcal{B}\left((d-\frac{\eta}{2},0), \frac{\eta}{4}\right)\right\} \subset A(\eta, d)$$

$$A_+(\eta, d) = \left\{ (z_1, z_2) \in \Delta_0(\overline{D})^2 : z_1 \in \mathcal{B}\left((0,0), \frac{\eta}{4}\right) \text{ and } z_2 \in \mathcal{B}\left((d+\frac{\eta}{2},0), \frac{\eta}{4}\right)\right\} \subset A(\eta, d).$$

For any $i \in \{2, \ldots, p^{m_1,m_2}\}$, when $\eta$ is small enough, the couple of points $(z_1, z_2) \in A(\eta, D_i^{m_1,m_2})$ (resp. $A_-(\eta, D_i^{m_1,m_2})$ and $A_+(\eta, D_i^{m_1,m_2})$) are such that $D_{i-1}^{m_1,m_2} < D_i^{m_1,m_2} - \eta < d(z_1,z_2) < D_i^{m_1,m_2} + \eta < D_{i+1}^{m_1,m_2}$ (resp. $D_{i-1}^{m_1,m_2} < D_i^{m_1,m_2} - \eta < d(z_1,z_2) < D_i^{m_1,m_2}$ and $D_i^{m_1,m_2} < d(z_1,z_2) < D_i^{m_1,m_2} + \eta < D_{i+1}^{m_1,m_2}$).

We now derive the following events

$$A_0 := \left\{ \varphi \in \overline{\Omega} : \varphi_{\Delta_0(\overline{D})} = \emptyset \right\}$$

$$A_i^{m_1,m_2}(\eta) = \left\{ \varphi \in \overline{\Omega} : \varphi_{\Delta_0(\overline{D})} = \{z_1^{m_1}, z_2^{m_2}\} \text{ with } (z_1,z_2) \in A(\eta, D_i^{m_1,m_2}) \right\}$$

$$A_{i,-}^{m_1,m_2}(\eta) = \left\{ \varphi \in \overline{\Omega} : \varphi_{\Delta_0(\overline{D})} = \{z_1^{m_1}, z_2^{m_2}\} \text{ with } (z_1,z_2) \in A_-(\eta, D_i^{m_1,m_2}) \right\} \subset A_i^{m_1,m_2}(\eta)$$

$$A_{i,+}^{m_1,m_2}(\eta) = \left\{ \varphi \in \overline{\Omega} : \varphi_{\Delta_0(\overline{D})} = \{z_1^{m_1}, z_2^{m_2}\} \text{ with } (z_1,z_2) \in A_+(\eta, D_i^{m_1,m_2}) \right\} \subset A_i^{m_1,m_2}(\eta).$$

First of all note that for any $\varphi_0 \in A_0$,

$$\boldsymbol{y}^T \boldsymbol{LPL}_{\overline{\Lambda}}^{(1)}(\varphi_0; \boldsymbol{\theta}^\star) = \sum_{m=1}^{M} y_1^{m,m} e^{-\theta_1^{\star,m,m}} |\overline{\Lambda}| - 0$$

For any $\varphi_{i,-}^{m_1,m_2} \in A_{i,-}^{m_1,m_2}$ and any $\varphi_{i,+}^{m_1,m_2} \in A_{i,+}^{m_1,m_2}$ for $i = 2, \ldots, p^{m_1,m_2}$

$$\boldsymbol{y}^T \boldsymbol{R}(\varphi_{i,-}^{m_1,m_2}) = y_1^{m_1,m_1} + y_1^{m_2,m_2} + 2y_i^{m_1,m_2}$$

$$\boldsymbol{y}^T \boldsymbol{R}(\varphi_{i,+}^{m_1,m_2}) = y_1^{m_1,m_1} + y_1^{m_2,m_2} + 2y_i^{m_1,m_2}(1 - \delta_{i,p^{m_1,m_2}})$$

We leave the reader to check that for every $\varepsilon > 0$ there exists $\eta > 0$ small enough such that

$$\left| \boldsymbol{y}^T \boldsymbol{L}(\varphi_{i,+}^{m_1,m_2}; \boldsymbol{\theta}^\star) - \boldsymbol{y}^T \boldsymbol{L}(\varphi_{i,-}^{m_1,m_2}; \boldsymbol{\theta}^\star) \right| \leq \varepsilon.$$

Therefore for every $\varepsilon > 0$ we have for $\eta$ small enough

$$\begin{aligned} 2|y^{m_1,m_2}_{p^{m_1,m_2}}| &= \Big| \boldsymbol{y}^T \left( \boldsymbol{L}(\varphi^{m_1,m_2}_{p^{m_1,m_2},+};\boldsymbol{\theta}^\star) - \boldsymbol{R}(\varphi^{m_1,m_2}_{p^{m_1,m_2},+}) \right) \\ &\quad - \boldsymbol{y}^T \left( \boldsymbol{L}(\varphi^{m_1,m_2}_{p^{m_1,m_2},-};\boldsymbol{\theta}^\star) - \boldsymbol{R}(\varphi^{m_1,m_2}_{p^{m_1,m_2},-}) \right) + 2y^{m_1,m_2}_{p^{m_1,m_2}} \Big| \le \varepsilon. \end{aligned}$$

By choosing $\varepsilon = |y^{m_1,m_2}_{p^{m_1,m_2}}|$, this leads to $y^{m_1,m_2}_{p^{m_1,m_2}} = 0$. By iterating this argument, we obtain that for any $m_1, m_2 \in \{1,\dots,M\}$, $y_2^{m_1,m_2} = \dots = y^{m_1,m_2}_{p^{m_1,m_2}} = 0$. It remains to prove that $y_1^{1,1} = \dots = y_1^{M,M} = 0$. For this, consider the following configuration set indexed by $n \ge 1$

$$A_n^{m_1} = \left\{ \varphi \in \overline{\Omega} : \varphi_{\Delta_0(\overline{D})} = \{z_1^{m_1},\dots,z_n^{m_1}\} \text{ with } z_1,\dots,z_n \in \Delta_0(\overline{D}) \right\}.$$

For any $\varphi_n^{m_1} \in A_n^{m_1}$, we have

$$\begin{aligned} \boldsymbol{y}^T \boldsymbol{L}(\varphi_n^{m_1,m_1};\boldsymbol{\theta}^\star) &= \sum_{m_1=1}^M y_1^{m_1,m_1} \int_{\overline{\Lambda}\times\mathbb{M}} e^{-V\left(x^m|\varphi_n^{m_1};\boldsymbol{\theta}^\star\right)} \mu(dx^m) - |\overline{\Lambda}| e^{-\theta_1^{\star,m_1,m_1}} \\ \boldsymbol{y}^T \boldsymbol{R}(\varphi_n^{m_1,m_1}) &= n y_1^{m_1,m_1} \end{aligned}$$

Hence for every $\varepsilon > 0$ we have for $n$ large enough by using the local stability property

$$\begin{aligned} |y_1^{m_1,m_1}| &= \left| \frac{1}{n} \boldsymbol{y}^T \left( \boldsymbol{L}(\varphi_n^{m_1,m_2};\boldsymbol{\theta}^\star) + \boldsymbol{R}(\varphi_n^{m_1,m_2}) \right) + y_1^{m_1,m_1} \right| \\ &\le \frac{2|\overline{\Lambda}|}{n} e^K \sum_{m_1=1}^M |y_1^{m_1,m_1}| \le \varepsilon. \end{aligned}$$

By choosing $\varepsilon = |y_1^{m_1,m_1}|/2$, this leads to $y_1^{m_1,m_1} = 0$.

### 5.3. *k−nearest-neighbour multi-Strauss marked point process*

Assumption **[Ident]** (resp. **[SDP]**) is proven without any change in the proof of the multi-Strauss marked point process for every $k \ge 1$ (resp. $k \ge 2$). The proof of **[SDP]** for $k = 1$ is omitted.

### 5.4. *Strauss disc process*

*Assumption [**Ident**]*

Consider

$$\begin{aligned} A_1 &:= \left\{ \varphi \in \Omega : \varphi_{\mathcal{B}(0,D)} = \emptyset \right\} \\ A_2 &:= \left\{ \varphi \in \Omega : \varphi_{\mathcal{B}(0,D)} = \{z^{m'}\},\ z \in \mathcal{B}(0,D/2),\ m' \in [D/2,D] \right\} \end{aligned}$$

and define $B_1 := \mathbb{M} \times A_1$ and $B_2 := \mathbb{M} \times A_2$. Then, for any $(m_1, \varphi_1, m_2, \varphi_2) \in B_1 \times B_2$

$$\underline{\boldsymbol{V}} = \begin{pmatrix} 1 & 0 \\ 1 & 1 \end{pmatrix},$$

which is injective since $\det(\underline{\boldsymbol{V}}) = 1$.

*Assumption* ***[SDP]***

Consider for $n \geq 1$

$$\begin{aligned} A_0 &:= \left\{ \varphi \in \overline{\Omega} : \varphi_{\Delta_0(\overline{D})} = \emptyset \right\} \\ A_n &:= \left\{ \varphi \in \overline{\Omega} : \varphi_{\Delta_0(\overline{D})} = \{ z_1^{m_1}, \ldots, z_n^{m_n} \}, \right. \\ &\qquad \left. z_1, \ldots, z_n \in \mathcal{B}(0, D/2),\ m_1, \ldots, m_n \in [D/2, D] \right\}. \end{aligned}$$

Note that for every $\varphi_0 \in A_0$ and any $\varphi_n \in A_n$

$$\boldsymbol{y}^T \boldsymbol{R}(\varphi_n) = n y_1 + n(n-1) y_2$$

Note also that from the local stability property $|\boldsymbol{y}^T \boldsymbol{L}(\varphi_n; \boldsymbol{\theta}^\star) \leq 2(|y_1| + |y_2|)|\overline{\Lambda}|e^K$. Then, for every $\varepsilon > 0$ we have for $n$ large enough

$$|y_2| = \left| \frac{1}{n(n-1)} \boldsymbol{y}^T \left( \boldsymbol{L}(\varphi_n; \boldsymbol{\theta}^\star) - \boldsymbol{R}_n(\varphi_n) \right) + y_2 \right| \leq \varepsilon.$$

By choosing $\varepsilon = |y_2|/2$, this leads to $y_2 = 0$. Then, for every $\varepsilon' > 0$ we have for $n$ large enough

$$|y_1| = \left| \frac{1}{n} (y_1, 0)^T \left( \boldsymbol{L}(\varphi_n; \boldsymbol{\theta}^\star) - \boldsymbol{R}_n(\varphi_n) \right) + y_1 \right| \leq \varepsilon'.$$

By choosing $\varepsilon' = |y_1|/2$, this leads to $y_1 = 0$.

### 5.5. Geyer's triplet point process

*Assumption* ***[Ident]***

By considering

$$A_1 := \left\{ \varphi \in \Omega : \varphi_{\mathcal{B}(0,D)} = \emptyset \right\},$$

$$A_2 := \left\{ \varphi \in \Omega : \varphi_{\mathcal{B}(0,D)} = \{z\},\ z \in \mathcal{B}((0, D/2), D/4), \varphi_{\mathcal{B}(0,2D) \setminus \mathcal{B}(0,D)} = \emptyset \right\}$$

and

$$A_3 := \left\{ \varphi \in \Omega : \varphi_{\mathcal{B}(0,D)} = \{z_1, z_2\}, z_1, z_2 \in \mathcal{B}((0, D/2), D/4), \varphi_{\mathcal{B}(0,2D) \setminus \mathcal{B}(0,D)} = \emptyset \right\},$$

we have for any $(\varphi_1, \varphi_2, \varphi_3) \in A_1 \times A_2 \times A_3$

$$\underline{\boldsymbol{V}} = \begin{pmatrix} 1 & 0 & 0 \\ 1 & 1 & 1 \\ 1 & 2 & 3 \end{pmatrix},$$

which is clearly injective since $\det \underline{\boldsymbol{V}} = 1$.

*Assumption* ***[SDP]***

Denote by $A_0$ any configuration set. Consider $A_n(\eta)$ for $n \geq 1$ and for some $0 < \eta < D$ the following configuration set

$$A_n(\eta) = \left\{ \varphi \in \overline{\Omega} : \varphi_{\Delta_0(\overline{D})} = \{z_1, \ldots, z_n\} \text{ with } z_1, \ldots, z_n \in \mathcal{B}(0, \eta) \right\}.$$

We leave the reader to chek that for $j = 1, \ldots, p$ and for any $\varphi_n \in A_n(\eta)$

$$\left| \int_{\overline{\Lambda}} v_j(x|\varphi_n) e^{-\boldsymbol{\theta}^{\star T} \boldsymbol{v}(x|\varphi_n)} dx \right| \leq n^{j-1} |\overline{\Lambda}| e^K,$$

where $K$ comes from the local stability property. Let us also remark that

$$\boldsymbol{y}^T \boldsymbol{R}(\varphi_n) = n y_1 + n(n-1) y_2 + \frac{n(n-1)(n-2)}{2} y_3 - \boldsymbol{y}^T \sum_{x \in \varphi_0 \cap \overline{\Lambda}} \boldsymbol{v}\left(x|\varphi_0 \setminus x\right).$$

Therefore by combining these two arguments, for every $\varepsilon > 0$, we have for $n$ large enough

$$\frac{|y_3|}{2} = \left| \frac{y_3}{2} + \frac{\boldsymbol{y}^T (\boldsymbol{L}(\varphi_n; \boldsymbol{\theta}^\star) - \boldsymbol{R}(\varphi_n))}{n(n-1)(n-2)} \right| \leq \varepsilon.$$

By choosing $\varepsilon = \frac{|y_3|}{4}$, this leads to $y_3 = 0$. Then, for every $\varepsilon'$, we may obtain for $n$ large enough

$$|y_2| = \left| y_2 + \frac{(y_1, y_2, 0)^T (\boldsymbol{L}(\varphi_n; \boldsymbol{\theta}^\star) - \boldsymbol{R}(\varphi_n))}{n(n-1)} \right| \leq \varepsilon'.$$

By choosing $\varepsilon' = |y_2|/2$, this leads to $y_2 = 0$. And then for every $\varepsilon''$ for $n$ large enough

$$|y_1| = \left| y_1 + \frac{(y_1, 0, 0)^T (\boldsymbol{L}(\varphi_n; \boldsymbol{\theta}^\star) - \boldsymbol{R}(\varphi_n))}{n} \right| \leq \varepsilon''.$$

By choosing $\varepsilon'' = \frac{|y_1|}{2}$, this finally leads to $y_1 = 0$.

### 5.6. Area-interaction model

*Assumption [**Ident**]*

By considering

$$A_1 := \left\{\varphi \in \Omega : \varphi_{\mathcal{B}(0,D)} = \emptyset\right\}$$

and for some small $\eta > 0$

$$A_2 := A_2(\eta) = \left\{\varphi \in \Omega : \varphi_{\mathcal{B}(0,D)} = \{z\}, z \in \mathcal{B}(0,\eta)\right\}.$$

we obtain for any $(\varphi_1, \varphi_2) \in A_1 \times A_2$

$$\underline{\boldsymbol{V}} = \begin{pmatrix} 1 & \pi R^2 \\ 1 & v_2(0|\varphi_2) \end{pmatrix}.$$

Since for any $\eta < R$, $0 < v_2(0|\varphi_2) < \pi R^2$, $\det(\underline{\boldsymbol{V}}) \neq 0$.

*Assumption [**SDP**]*

Consider $A_k(\eta)$ for $k = 0, 1, 2$ and for some $0 < \eta < D$ the following configuration set

$$A_k(\eta) = \left\{\varphi \in \overline{\Omega} : \varphi_{\Delta_0(\overline{D})} = \{z_1, \ldots, z_{k+1}\} \text{ with } z_1, \ldots, z_{k+1} \in \mathcal{B}(0,\eta)\right\}$$

For any $\varphi_k \in A_k(\eta)$ for $k = 0, 1, 2$,

$$\begin{aligned} R(\varphi_1) &= 2y_1 + y_2 g_1(\varphi_1)) - (y_1 + y_2 \pi R^2) \\ R(\varphi_2) &= 3y_1 + y_2 g_2(\varphi_2) - (y_1 + y_2 \pi R^2) \end{aligned}$$

For every $\varepsilon > 0$, there exists $\eta > 0$ such that for any $\varphi_k \in A_k(\eta)$ $(k = 1, 2)$, $|\boldsymbol{y}^T \boldsymbol{L}(\varphi_k; \boldsymbol{\theta}^\star)| \leq \varepsilon$ and such that $|y_2 g_k(\varphi_k)| \leq \varepsilon$. Then,

$$|y_1| = \left|y_1 - \boldsymbol{y}^T(L(\varphi_1; \boldsymbol{\theta}^\star) - L_2(\varphi_2; \boldsymbol{\theta}^\star)) - (R(\varphi_1) - R(\varphi_2))\right| \leq 4\varepsilon.$$

By choosing $\varepsilon = |y_1|/8$, this leads to $y_1 = 0$. Now,

$$|y_2|\pi R^2 = \left|y_2 \pi R^2 - (0, y_2)^T (L(\varphi_1; \boldsymbol{\theta}^\star) - R(\varphi_1))\right| \leq 2\varepsilon.$$

And by choosing $\varepsilon = \pi R^2 |y_2|/4$, this finally leads to $y_2 = 0$.

## 6. Annex: Proofs

### 6.1. Tools

Let us start by presenting a particular case of the Campbell Theorem combined with the Glötzl Theorem that is widely used in our future proofs.

**Corollary 3.** *Assume that the (marked) point process $\Phi$ with probability measure $P$ is stationary. Let $\Lambda$ be a bounded Borel set, let $\varphi \in \Omega$ and let $g$ be a function satisfying $g(x^m, \varphi_x) = g(0^m, \varphi)$ for all $x^m \in \$$. Define $M$ a random variable with its distribution $\lambda^{\mathbb{m}}$ and $f(m, \varphi) = g(0^m, \varphi)e^{-V(0^m|\varphi)}$ and assume that $f \in L^1(\lambda^{\mathbb{m}} \otimes P)$. Then,*

$$\boldsymbol{E}\Big( \sum_{x^m \in \Phi_\Lambda} g(x^m, \Phi \setminus x^m) \Big) = |\Lambda| \, \boldsymbol{E}\Big( g\left(0^M, \Phi\right) e^{-V\left(0^M|\Phi\right)} \Big) \tag{6.1}$$

**Proof.**

$$\begin{aligned}
\boldsymbol{E}\Big( \sum_{x^m \in \Phi_\Lambda} g(x^m, \Phi \setminus x^m) \Big) &= \int_\Omega \sum_{x^m \in \varphi} g(x^m, \varphi \setminus x^m) \mathbb{1}_\Lambda(x) P(d\varphi) \\
&= \int_{\$} \int_\Omega g(x^m, \varphi) \mathbb{1}_\Lambda(x) P^!_{x^m}(d\varphi) N_P(dx^m) \\
&= \int_{\$} \int_\Omega g(x^m, \varphi) \mathbb{1}_\Lambda(x) \nu_P(x^m) P^!_{x^m}(d\varphi) \mu(dx^m) \\
&= \int_{\$} \int_\Omega g(x^m, \varphi) \mathbb{1}_\Lambda(x) e^{-V(x^m|\varphi)} P(d\varphi) \mu(dx^m) \\
&= \int_{\Lambda \times \mathbb{M}} \int_\Omega g(x^m, \varphi) e^{-V(x^m|\varphi)} P(d\varphi) \mu(dx^m) \\
&= |\Lambda| \int_{\mathbb{M}} \int_\Omega g(x^m, \varphi) e^{-V(x^m|\varphi)} P(d\varphi) \lambda^{\mathbb{m}}(dm) \\
&= |\Lambda| \boldsymbol{E}\Big( g(0^M, \Phi) e^{-V\left(0^M|\Phi\right)} \Big)
\end{aligned}$$

where $\nu_P(\cdot)$ is the Radon-Nikodym derivative of $N_P$ with respect to $\mu$. ■

Let us now present a version of an ergodic theorem obtained by Nguyen and Zessin (1979) and widely used in this paper. Let $\widetilde{D} > 0$ and denote by $\Delta_0$ the following fixed domain

**Theorem 4** (Nguyen and Zessin (1979))**.** *Let $\{H_G, G \in \mathcal{B}_b^2\}$ be a family of random variables, which is covariant, that for all $x \in \mathbb{R}^2$,*

$$H_{G_x}(\varphi_x) = H_G(\varphi), \quad \text{for a.e. } \varphi$$

*and additive, that is for every disjoint $G_1, G_2 \in \mathcal{B}_b^2$,*

$$H_{G_1 \cup G_2} = H_{G_1} + H_{G_2}, \quad a.s.$$

*Let $\mathcal{I}$ be the sub$-\sigma-$algebra of $\mathcal{F}$ consisting of translation invariant (with probability 1) sets. Assume there exists a nonnegative and integrable random variable $Y$ such that $|H_G| \leq Y$ a.s. for every convex $G \subset \Delta_0$. Then,*

$$\lim_{n \to +\infty} \frac{1}{|G_n|} H_{G_n} = \frac{1}{|\Delta_0|} E(H_{\Delta_0}|\mathcal{I}), \quad a.s.$$

*for each regular sequence $G_n \to \mathbb{R}^2$.*

### 6.2. Proof of Theorem 1

Due to the decomposition of stationary measures as a mixture of ergodic measures (see Preston (1976)), one only needs to prove Theorem 1 by assuming that $P_{\boldsymbol{\theta}^\star}$ is ergodic. From now on, $P_{\boldsymbol{\theta}^\star}$ is assumed to be ergodic. The tool used to obtain the almost sure convergence is a convergence theorem for minimum contrast estimators established by Guyon (1995).

We proceed in three stages.

*Step 1. Convergence of $U_n(\Phi;\boldsymbol{\theta})$.*

Decompose $U_n(\varphi;\boldsymbol{\theta}) = \frac{1}{|\Lambda_n|}\left(H_{1,\Lambda_n}(\varphi) + H_{2,\Lambda_n}(\varphi)\right)$ with

$$H_{1,\Lambda_n}(\varphi) = \int_{\Lambda_n\times\mathbb{M}} e^{-V(x^m|\varphi;\boldsymbol{\theta})}\mu(dx^m)$$

and

$$H_{2,\Lambda_n}(\varphi) = \sum_{x^m\in\Phi_{\Lambda_n}} V\left(x^m|\varphi\setminus x^m;\boldsymbol{\theta}\right).$$

Under the assumption **[Mod]**, one can apply Theorem 4 (Nguyen and Zessin (1979)) to the process $H_{1,\Lambda_n}$. And from Corollary 3, we obtain $P_{\boldsymbol{\theta}^\star}-$almost surely as $n\to+\infty$

$$\frac{1}{|\Lambda_n|}H_{1,\Lambda_n} \to \boldsymbol{E}\Big(\int_{\mathbb{M}} e^{-V(0^m|\Phi;\boldsymbol{\theta})}\lambda^{\mathrm{m}}(dm)\Big) = \boldsymbol{E}\Big(e^{-V\left(0^M|\Phi;\boldsymbol{\theta}\right)}\Big). \tag{6.2}$$

Now, let $G\subset\Delta_0$, we clearly have

$$|H_{2,G}(\varphi)| \le \sum_{x^m\in\Phi_G}|V\left(x^m|\varphi\setminus x^m;\boldsymbol{\theta}\right)| \le \sum_{x^m\in\varphi_{\Delta_0}}|V\left(x^m|\varphi\setminus x^m;\boldsymbol{\theta}\right)|$$

Under the assumption **[Mod]** and from Corollary 3, we have

$$\boldsymbol{E}\left(\sum_{x^m\in\Phi_{\Delta_0}}|V\left(x^m|\Phi\setminus x^m;\boldsymbol{\theta}\right)|\right) = |\Delta_0|\boldsymbol{E}\left(\left|V\left(0^M|\Phi;\boldsymbol{\theta}\right)\right|e^{-V\left(0^M|\Phi;\boldsymbol{\theta}^\star\right)}\right) < +\infty$$

This means that for all $G\subset\Delta_0$, there exists a random variable $Y\in L^1(P_{\boldsymbol{\theta}^\star})$ such that $|H_{2,G}|\le Y$. Thus, under the ssumption **[Mod]** and from Theorem 4 (Nguyen and Zessin (1979)) and from Corollary 3, we have $P_{\boldsymbol{\theta}^\star}-$almost surely

$$\begin{aligned}\frac{1}{|\Lambda_n|}H_{2,\Lambda_n} \quad &\to \quad \frac{1}{|\Delta_0|}\boldsymbol{E}\Big(\sum_{x^m\in\Phi_{\Delta_0}}V\left(x^m|\Phi\setminus x^m;\boldsymbol{\theta}\right)\Big) \\ &= \quad \boldsymbol{E}\left(V\left(0^M|\Phi;\boldsymbol{\theta}\right)e^{-V\left(0^M|\Phi;\boldsymbol{\theta}^\star\right)}\right).\end{aligned} \tag{6.3}$$

We have the result by combining (6.2) and (6.3): $P_{\boldsymbol{\theta}^\star}-$almost surely

$$U_n(\cdot;\boldsymbol{\theta}) \to U(\boldsymbol{\theta}) = \boldsymbol{E}\Big(e^{-V\left(0^M|\Phi;\boldsymbol{\theta}\right)} + V\left(0^M|\Phi;\boldsymbol{\theta}\right)e^{-V\left(0^M|\Phi;\boldsymbol{\theta}^\star\right)}\Big) \tag{6.4}$$

*Step 2. $U_n(\Phi;\cdot)$ a contrast function*

Recall that $U_n(\cdot;\boldsymbol{\theta})$ is a contrast function if there exists a function $K(\cdot,\boldsymbol{\theta}^\star)$ (i.e. nonnegative function equal to zero if and only if $\boldsymbol{\theta}=\boldsymbol{\theta}^\star$) such that $P_{\boldsymbol{\theta}^\star}-$almost surely $U_n(\varphi;\boldsymbol{\theta})-U_n(\varphi;\boldsymbol{\theta})\to K(\boldsymbol{\theta},\boldsymbol{\theta}^\star)$. From Step 1, we have as $n\to+\infty$

$$\begin{aligned} K(\boldsymbol{\theta},\boldsymbol{\theta}^\star) &= \boldsymbol{E}\Big(e^{-V\left(0^M|\Phi;\boldsymbol{\theta}^\star\right)}\Big(e^{V\left(0^M|\Phi;\boldsymbol{\theta}\right)-V\left(0^M|\Phi;\boldsymbol{\theta}^\star\right)} \\ &\quad -\Big(1+V\left(0^M|\Phi;\boldsymbol{\theta}\right)-V\left(0^M|\Phi;\boldsymbol{\theta}^\star\right)\Big)\Big)\Big) \\ &= \boldsymbol{E}\Big(e^{-\boldsymbol{\theta}^{\star T}\boldsymbol{v}\left(0^M|\varphi\right)}\Big(e^{(\boldsymbol{\theta}-\boldsymbol{\theta}^\star)^T\boldsymbol{v}\left(0^M|\Phi\right)}-(1+(\boldsymbol{\theta}-\boldsymbol{\theta}^\star)^T\boldsymbol{v}\left(0^M|\Phi\right))\Big)\Big). \end{aligned} \tag{6.5}$$

Let $\boldsymbol{y}\in\mathbb{R}^{p+1}\setminus\{0\}$, and assume $\boldsymbol{y}^T\boldsymbol{v}\left(0^m|\varphi\right)=0$ for $\lambda^{\mathrm{m}}\otimes P_{\boldsymbol{\theta}^\star}-$a.e. $(m,\varphi)$. By assuming **[Ident]**, it follows that for $i=1,\ldots,\ell$ $(\ell\geq p)$ $\boldsymbol{y}^T\boldsymbol{v}\left(0^{m_i}|\varphi_i\right)=0$ for all $(m_i,\varphi_i)\in B_i$. From the injectivity of the matrix with entries $v_j(0^{m_i}|\varphi_i)$ for all $(m_1,\varphi_1,\ldots,m_\ell,\varphi_\ell)\in B_1\times\ldots\times B_\ell$ it comes that $\boldsymbol{y}=0$ which leads to a contradiction. Therefore, for $\boldsymbol{\theta}\neq\boldsymbol{\theta}^\star$, the assertion $(\boldsymbol{\theta}-\boldsymbol{\theta}^\star)^T\boldsymbol{v}\left(0^m|\varphi\right)\neq 0$ holds for $\lambda^{\mathrm{m}}\otimes P_{\boldsymbol{\theta}^\star}-$a.e. $(m,\varphi)$.

By noticing that the function $t\mapsto e^t-(1+t)$ is nonnegative and is equal to zero if and only if $t=0$, one concludes that the random variable $e^{-\boldsymbol{\theta}^{\star T}\boldsymbol{v}\left(0^M|\Phi\right)}\times\Big(e^{(\boldsymbol{\theta}-\boldsymbol{\theta}^\star)^T\boldsymbol{v}\left(0^M|\Phi\right)}-(1+(\boldsymbol{\theta}-\boldsymbol{\theta}^\star)^T\boldsymbol{v}\left(0^M|\Phi\right))\Big)$ is almost surely positive for $\boldsymbol{\theta}\neq\boldsymbol{\theta}^\star$ and equals to zero when $\boldsymbol{\theta}=\boldsymbol{\theta}^\star$.

Before ending this step, note that for any $\varphi$, $U_n(\varphi;\cdot)$ and $K(\cdot,\boldsymbol{\theta}^\star)$ are clearly continuous functions.

*Step 3. Modulus of continuity.*

The modulus of continuity of the contrast process defined for all $\varphi\in\Omega$ and all $\eta>0$ by

$$W_n(\varphi,\eta)=\sup\Big\{\Big|U_n(\varphi;\boldsymbol{\theta})-U_n(\varphi;\boldsymbol{\theta}')\Big|:\boldsymbol{\theta},\boldsymbol{\theta}'\in\boldsymbol{\Theta},||\boldsymbol{\theta}-\boldsymbol{\theta}'||\leq\eta\Big\}$$

is such that there exists a sequence $(\varepsilon_k)_{k\geq1}$, with $\varepsilon_k\to0$ as $k\to+\infty$ such that for all $k\geq1$

$$P\left(\limsup_{n\to+\infty}\left(W_n\left(\Phi,\frac{1}{k}\right)\geq\varepsilon_k\right)\right)=0. \tag{6.6}$$

Let us start to write $W_n\left(\varphi,\frac{1}{k}\right)\leq W_{1,n}\left(\varphi,\frac{1}{k}\right)+W_{2,n}\left(\varphi,\frac{1}{k}\right)$ with

$$\begin{aligned} W_{1,n}\left(\varphi,\frac{1}{k}\right) &= \sup\Big\{\Big|\frac{1}{|\Lambda_n|}\int_{\Lambda_n\times\mathbb{M}}\Big(e^{-V(x^m|\varphi;\boldsymbol{\theta})}-e^{-V\left(x^m|\varphi;\boldsymbol{\theta}'\right)}\Big)\mu(dx^m)\Big| \\ &\qquad :\boldsymbol{\theta},\boldsymbol{\theta}'\in\boldsymbol{\Theta},||\boldsymbol{\theta}-\boldsymbol{\theta}'||\leq\frac{1}{k}\Big\} \end{aligned}$$

and

$$W_{2,n}\left(\varphi,\frac{1}{k}\right) = \sup\left\{\Big|\sum_{x^m\in\varphi_{\Lambda_n}} V\left(x^m|\varphi\setminus x^m;\boldsymbol{\theta}\right)-V\left(x^m|\varphi\setminus x^m;\boldsymbol{\theta}'\right)\Big| \right.$$
$$\left. :\boldsymbol{\theta},\boldsymbol{\theta}'\in\boldsymbol{\Theta},||\boldsymbol{\theta}-\boldsymbol{\theta}'||\leq\frac{1}{k}\right\}.$$

From the assumption **[Mod]** it comes

$$\begin{aligned} W_{1,n}\left(\varphi,\frac{1}{k}\right) &\leq \sup\left\{\frac{1}{|\Lambda_n|}\int_{\Lambda_n\times\mathbb{M}}\left(|(\boldsymbol{\theta}-\boldsymbol{\theta}')^T\boldsymbol{v}\left(x^m|\varphi\right)|e^K\right)\mu(dx^m)\right. \\ &\qquad\left. :\boldsymbol{\theta},\boldsymbol{\theta}'\in\boldsymbol{\Theta},||\boldsymbol{\theta}-\boldsymbol{\theta}'||\leq\frac{1}{k}\right\} \\ &\leq \frac{e^K}{k}\frac{1}{|\Lambda_n|}\int_{\Lambda_n\times\mathbb{M}}||\boldsymbol{v}\left(x^m|\Phi\right)||\mu(dx^m). \end{aligned}$$

Under the assumption **[Mod]**, one can apply Theorem 4 (Nguyen and Zessin (1979)) to obtain as $n\to+\infty$

$$\frac{1}{|\Lambda_n|}\int_{\Lambda_n\times\mathbb{M}}||\boldsymbol{v}\left(x^m|\Phi\right)||\mu(dx^m)\to\boldsymbol{E}\left(||\boldsymbol{v}\left(0^M|\Phi\right)||\right),\ \text{for } P_{\boldsymbol{\theta}^\star}-a.e.\,\varphi.$$

So there exists $n_0^{(1)}(k)$ such that for all $n\geq n_0^{(1)}(k)$ we have

$$W_{1,n}\left(\varphi,\frac{1}{k}\right)\leq\frac{2e^K}{k}\boldsymbol{E}\left(||\boldsymbol{v}\left(0^M|\Phi\right)||\right),$$

for $P_{\boldsymbol{\theta}^\star}$-a.e. $\varphi$. By using the same arguments, one can prove that there exists $n_0^{(2)}(k)$ such that for all $n\geq n_0^{(2)}(k)$ we have

$$W_{2,n}\left(\varphi,\frac{1}{k}\right)\leq\frac{2}{k}\boldsymbol{E}\left(||\boldsymbol{v}\left(0^M|\Phi\right)||e^{-V\left(0^M|\Phi;\boldsymbol{\theta}^\star\right)}\right)\leq\frac{2e^K}{k}\boldsymbol{E}\left(||\boldsymbol{v}\left(0^M|\Phi\right)||\right),$$

for $P_{\boldsymbol{\theta}^\star}-$a.e. $\varphi$. And so for all $n\geq n_0(k)=\max(n_0^{(1)}(k),n_0^{(2)}(k))$, we have $P_{\boldsymbol{\theta}^\star}-$a.s.

$$W_n\left(\varphi,\frac{1}{k}\right)\leq W_{1,n}\left(\varphi,\frac{1}{k}\right)+W_{2,n}\left(\varphi,\frac{1}{k}\right)<\frac{\delta}{k},\ \text{for } P_{\boldsymbol{\theta}^\star}-a.e.\,\varphi.$$

with $\delta=4e^K\boldsymbol{E}\left(||\boldsymbol{v}\left(0^M|\Phi\right)||\right)$. Since,

$$\begin{aligned} \limsup_{n\to+\infty}\left\{W_n\left(\varphi,\frac{1}{k}\right)\geq\frac{\delta}{k}\right\} &= \bigcap_{m\in\mathbb{N}}\bigcup_{n\geq m}\left\{W_n\left(\varphi,\frac{1}{k}\right)\geq\frac{\delta}{k}\right\} \\ &\subset \bigcup_{n\geq n_0(k)}\left\{W_n\left(\varphi,\frac{1}{k}\right)\geq\frac{\delta}{k}\right\},\ \text{for } P_{\boldsymbol{\theta}^\star}-a.e.\,\varphi. \end{aligned}$$

the expected result (6.6) is proved.
*Conclusion step.* The Steps 1, 2 and 3 ensure the fact that we can apply a consistency result on minimum contrast estimators, see Guyon (1995).

### 6.3. Proof of Theorem 2

*Step 1. Asymptotic normality of* $\boldsymbol{U}_n^{(1)}(\Phi;\boldsymbol{\theta}^\star)$

The aim is to prove that for any fixed $\widetilde{D}$, the following convergence in distribution as $n \to +\infty$

$$|\Lambda_n|^{1/2}\, \boldsymbol{U}_n^{(1)}(\Phi;\boldsymbol{\theta}^\star) \to \mathcal{N}\left(0, \underline{\boldsymbol{\Sigma}}(\widetilde{D},\boldsymbol{\theta}^\star)\right) \tag{6.7}$$

where the matrix $\underline{\boldsymbol{\Sigma}}(\widetilde{D},\boldsymbol{\theta}^\star)$ is defined by (4.4).

The idea is to apply to $\boldsymbol{U}_n^{(1)}(\Phi;\boldsymbol{\theta}^\star)$ a central limit theorem obtained by Jensen and Künsch (1994), Theorem 2.1. The following conditions have to be fulfilled to apply this result. For all $j = 1,\dots,p+1$

$(i)$ For all $i \in \mathbb{Z}^2$, $\boldsymbol{E}\left(\left(\boldsymbol{LPL}_{\Delta_i}^{(1)}(\Phi;\boldsymbol{\theta}^\star)\right)_j |\Phi_{\Delta_i^c}\right) = 0.$

$(ii)$ For all $i \in \mathbb{Z}^2$, $\boldsymbol{E}\left(\left|\left(\boldsymbol{LPL}_{\Delta_i}^{(1)}(\Phi;\boldsymbol{\theta}^\star)\right)_j\right|^3\right) < +\infty.$

$(iii)$ The matrix $\mathbb{V}\mathrm{ar}\left(|\Lambda_n|^{1/2}\boldsymbol{U}_n^{(1)}(\Phi;\boldsymbol{\theta}^\star)\right)$ converges to the matrix $\underline{\boldsymbol{\Sigma}}(\widetilde{D},\boldsymbol{\theta}^\star)$.

<u>Condition $(i)$</u>: From the stationarity of the process, it is sufficient to prove that

$$\boldsymbol{E}\left(\left(\boldsymbol{LPL}_{\Delta_0}^{(1)}(\Phi;\boldsymbol{\theta}^\star)\right)_j |\Phi_{\Delta_0^c}\right) = 0.$$

Recall that for any configuration $\varphi$

$$\begin{aligned}\left(\boldsymbol{LPL}_{\Delta_0}^{(1)}(\varphi;\boldsymbol{\theta}^\star)\right)_j &= -\int_{\Delta_0\times\mathbb{M}} v_j(x^m|\varphi)e^{-V(x^m|\varphi;\boldsymbol{\theta}^\star)}\mu(dx^m) \\ &\quad + \int_{\Delta_0\times\mathbb{M}} v_j(x^m|\varphi\setminus x^m)\varphi(dx^m).\end{aligned} \tag{6.8}$$

Denote respectively by $G_1(\varphi)$ and $G_2(\varphi)$ the first and the second right-hand term of (6.8) and by $E_i = \boldsymbol{E}\left(G_i(\Phi)|\Phi_{\Delta_0^c} = \varphi_{\Delta_0^c}\right)$. From the definition of Gibbs point processes,

$$E_2 = \frac{1}{Z_{\Delta_0}(\varphi_{\Delta_0^c})}\int_{\Omega_{\Delta_0}} Q(d\varphi_{\Delta_0})\int_{\mathbb{S}} \varphi_{\Delta_0}(dx^m)\mathbf{1}_{\Delta_0}(x)v_j(x^m|\varphi\setminus x^m)e^{-V\left(\varphi_{\Delta_0}|\varphi_{\Delta_0^c};\boldsymbol{\theta}^\star\right)}.$$

Denote by $\varphi' = (\varphi_{\Delta_0}, \varphi'_{\Delta_0^c})$. Since $Q$ is a Poisson process we can write

$$\begin{aligned}E_2 &= \frac{1}{Z_{\Delta_0}(\varphi_{\Delta_0^c})}\int_{\Omega} Q(d\varphi')\int_{\mathbb{S}} \varphi'(dx^m)\mathbf{1}_{\Delta_0}(x)v_j(x^m|\varphi\setminus x^m)e^{-V\left(\varphi_{\Delta_0}|\varphi_{\Delta_0^c};\boldsymbol{\theta}^\star\right)} \\ &= \frac{1}{Z_{\Delta_0}(\varphi_{\Delta_0^c})}\int_{\Omega} Q(d\varphi') \\ &\quad\times \int_{\mathbb{S}} \varphi'(dx^m)\mathbf{1}_{\Delta_0}(x)v_j(x^m|\varphi'_{\Delta_0}\cup\varphi_{\Delta_0^c}\setminus x^m)e^{-V\left(\varphi'_{\Delta_0}|\varphi_{\Delta_0^c};\boldsymbol{\theta}^\star\right)}\end{aligned}$$

Now, from Campbell Theorem (applied to the Poisson measure $Q$)

$$E_2 = \frac{1}{Z_{\Delta_0}(\varphi_{\Delta_0^c})} \int_{\Delta_0 \times \mathbb{M}} \mu(dx^m) \times \int_{\Omega} Q^!_{x^m}(d\varphi') v_j(x^m | \varphi'_{\Delta_0} \cup \varphi_{\Delta_0^c}) e^{-V\left(\varphi'_{\Delta_0} \cup x^m | \varphi_{\Delta_0^c}; \boldsymbol{\theta}^\star\right)}.$$

Since from Slivnyak-Mecke Theorem, $Q = Q^!_x$, one can obtain

$$\begin{aligned} E_2 &= \frac{1}{Z_{\Delta_0}(\varphi_{\Delta_0^c})} \int_{\Omega} Q(d\varphi') \int_{\Delta_0 \times \mathbb{M}} \mu(dx^m)\, v_j(x^m | \varphi'_{\Delta_0} \cup \varphi_{\Delta_0^c}) e^{-V\left(\varphi'_{\Delta_0} \cup x^m | \varphi_{\Delta_0^c}; \boldsymbol{\theta}^\star\right)} \\ &= \frac{1}{Z_{\Delta_0}(\varphi_{\Delta_0^c})} \int_{\Omega} Q(d\varphi_{\Delta_0}) \int_{\Delta_0 \times \mathbb{M}} \mu(dx^m) v_j(x^m | \varphi) e^{-V(x^m | \varphi; \boldsymbol{\theta}^\star)} e^{-V\left(\varphi_{\Delta_0} | \varphi_{\Delta_0^c}; \boldsymbol{\theta}^\star\right)} \\ &= -E_1 \end{aligned}$$

<u>Condition $(ii)$</u>: For any bounded domain $\Delta$ and any configuration $\varphi$, one may write for $j = 1, \dots, p+1$

$$\left| \left( \boldsymbol{LPL}^{(1)}_{\Delta}(\varphi; \boldsymbol{\theta}^\star) \right)_j \right|^3 \leq 4 \left| \int_{\Delta \times \mathbb{M}} v_j(x^m | \varphi) e^{-V(x^m | \varphi; \boldsymbol{\theta}^\star)} \mu(dx^m) \right|^3 + 4 \left| \sum_{x^m \in \varphi_\Delta} v_j(x^m | \varphi \setminus x^m) \right|^3$$

From **[Mod]**, both right-hand terms are integrable with respect to $P_{\boldsymbol{\theta}^\star}$, which implies that for any domain $\Delta$ and in particular for $\Delta_i$,

$$\boldsymbol{E}\left( \left| \left( \boldsymbol{LPL}^{(1)}_{\Delta_i}(\Phi; \boldsymbol{\theta}^\star) \right)_j \right|^3 \right) < +\infty.$$

<u>Condition $(iii)$</u>: let us start by noting that the vector $\boldsymbol{LPL}^{(1)}_{\Delta_i}(\varphi; \boldsymbol{\theta}^\star)$ depends only on $\varphi_{\Delta_j}$ for $j \in \mathbb{B}\left(i, \left\lceil \frac{D}{\tilde{D}} \right\rceil\right)$. Let

$$E_{i,j} := \boldsymbol{E}\left( \boldsymbol{LPL}^{(1)}_{\Delta_i}(\Phi; \boldsymbol{\theta}^\star)\, \boldsymbol{LPL}^{(1)}_{\Delta_j}(\Phi; \boldsymbol{\theta}^\star)^T \right) = E_{0,j-i},$$

by using the stationarity of the process. From definitions, we can obtain

$$\begin{aligned} &\mathbb{V}\mathrm{ar}\left( |\Lambda_n|^{1/2} \boldsymbol{U}^{(1)}_n(\Phi; \boldsymbol{\theta}^\star) \right) \\ &= |\Lambda_n|^{-1} \mathbb{V}\mathrm{ar}\left( \sum_{i \in I_n} \boldsymbol{LPL}^{(1)}_{\Delta_i}(\Phi; \boldsymbol{\theta}^\star) \right) \\ &= |\Lambda_n|^{-1} \sum_{i,j \in I_n} E_{i,j} \\ &= |\Lambda_n|^{-1} \sum_{i \in I_n} \left( \sum_{j \in I_n \cap \mathbb{B}\left(i, \left\lceil \frac{D}{\tilde{D}} \right\rceil\right)} E_{i,j} + \sum_{j \in I_n \cap \mathbb{B}\left(i, \left\lceil \frac{D}{\tilde{D}} \right\rceil\right)^c} E_{i,j} \right). \end{aligned}$$

By using condition $(i)$, one may assert that for any $j \in I_n \cap \mathbb{B}\left(i, \left\lceil \frac{D}{\widetilde{D}} \right\rceil\right)^c$

$$\begin{aligned}
E &:= \boldsymbol{E}\left(\boldsymbol{LPL}^{(1)}_{\Delta_i}(\Phi;\boldsymbol{\theta}^\star)\,\boldsymbol{LPL}^{(1)}_{\Delta_j}(\Phi;\boldsymbol{\theta}^\star)^T\right)\\
&= \boldsymbol{E}\left(\boldsymbol{E}\left(\boldsymbol{LPL}^{(1)}_{\Delta_i}(\Phi;\boldsymbol{\theta}^\star)\,\boldsymbol{LPL}^{(1)}_{\Delta_j}(\Phi;\boldsymbol{\theta}^\star)^T|\Phi_{\Delta_j}\right)\right)\\
&= \boldsymbol{E}\left(\boldsymbol{E}\left(\boldsymbol{LPL}^{(1)}_{\Delta_i}(\Phi;\boldsymbol{\theta}^\star)|\Phi_{\Delta_j}\right)\boldsymbol{LPL}^{(1)}_{\Delta_j}(\Phi;\boldsymbol{\theta}^\star)^T\right)\\
&= 0
\end{aligned}$$

Denote by $\widetilde{I}_n$ the following set

$$\widetilde{I}_n = I_n \cap \left(\cup_{i\in\partial I_n}\mathbb{B}\left(i, \left\lceil \frac{D}{\widetilde{D}} \right\rceil\right)\right).$$

We now obtain

$$\begin{aligned}
\mathbb{V}\mathrm{ar}\left(|\Lambda_n|^{1/2}\boldsymbol{U}^{(1)}_n(\Phi;\boldsymbol{\theta}^\star)\right) &= |\Lambda_n|^{-1}\sum_{i\in I_n}\sum_{j\in I_n\cap\mathbb{B}\left(i,\left\lceil\frac{D}{\widetilde{D}}\right\rceil\right)} E_{i,j}\\
&= |\Lambda_n|^{-1}\Bigg(\sum_{i\in I_n\setminus\widetilde{I}_n}\sum_{j\in I_n\cap\mathbb{B}\left(i,\left\lceil\frac{D}{\widetilde{D}}\right\rceil\right)} E_{i,j}\\
&\qquad + \sum_{i\in\widetilde{I}_n}\sum_{j\in I_n\cap\mathbb{B}\left(i,\left\lceil\frac{D}{\widetilde{D}}\right\rceil\right)} E_{i,j}\Bigg)
\end{aligned}$$

Using the stationarity and the definition of the domain $\Lambda_n$, one obtains

$$\begin{aligned}
|\Lambda_n|^{-1}\sum_{i\in I_n\setminus\widetilde{I}_n}\sum_{j\in I_n\cap\mathbb{B}\left(i,\left\lceil\frac{D}{\widetilde{D}}\right\rceil\right)} E_{i,j} &= |\Lambda_n|^{-1}|I_n\setminus\widetilde{I}_n|\sum_{j\in\mathbb{B}\left(0,\left\lceil\frac{D}{\widetilde{D}}\right\rceil\right)} E_{0,j}\\
&\to \underline{\boldsymbol{\Sigma}}(\widetilde{D},\boldsymbol{\theta}^\star),
\end{aligned}$$

as $n\to+\infty$ and

$$|\Lambda_n|^{-1}\left|\sum_{i\in\widetilde{I}_n}\sum_{j\in I_n\cap\mathbb{B}\left(i,\left\lceil\frac{D}{\widetilde{D}}\right\rceil\right)} E_{i,j}\right| \le |\Lambda_n|^{-1}|\widetilde{I}_n|\sum_{j\in\mathbb{B}\left(0,\left\lceil\frac{D}{\widetilde{D}}\right\rceil\right)}|E_{0,j}| \to 0$$

as $n \to +\infty$. Hence as $n \to +\infty$

$$\begin{aligned} \mathbb{V}\text{ar}\left(|\Lambda_n|^{1/2}\boldsymbol{U}_n^{(1)}(\Phi;\boldsymbol{\theta}^\star)\right) &= |\Lambda_n|^{-1}\sum_{i\in I_n}\sum_{j\in I_n\cap\mathbb{B}\left(i,\left\lceil\frac{\widetilde{D}}{D}\right\rceil\right)} E_{i,j} \\ &\sim \underbrace{|I_n||\Lambda_n|^{-1}}_{\widetilde{D}^{-2}}\sum_{k\in\mathbb{B}\left(0,\left\lceil\frac{\widetilde{D}}{D}\right\rceil\right)} E_{0,k} = \underline{\boldsymbol{\Sigma}}(\widetilde{D},\boldsymbol{\theta}^\star). \end{aligned} \quad (6.9)$$

*Step 2. Domination of $\underline{\boldsymbol{U}}_n^{(2)}(\Phi;\boldsymbol{\theta})$ in a neighborhood of $\boldsymbol{\theta}^\star$ and convergence of $\underline{\boldsymbol{U}}_n^{(2)}(\Phi;\boldsymbol{\theta}^\star)$*

Let $j,k=1,\dots,p$, recall that $\left(\underline{\boldsymbol{U}}_n^{(2)}(\varphi;\boldsymbol{\theta})\right)_{j,k}$ is defined for any configuration $\varphi$ by

$$\left(\underline{\boldsymbol{U}}_n^{(2)}(\varphi;\boldsymbol{\theta})\right)_{j,k} = \frac{1}{|\Lambda_n|}\int_{\Lambda_n\times\mathbb{M}} v_j(x^m|\varphi)v_k(x^m|\varphi)e^{-V(x^m|\varphi;\boldsymbol{\theta})}\mu(dx^m).$$

Using the local stability property it comes

$$\left(\underline{\boldsymbol{U}}_n^{(2)}(\varphi;\boldsymbol{\theta})\right)_{j,k} \le \frac{e^K}{|\Lambda_n|}\int_{\Lambda_n\times\mathbb{M}} v_j(x^m|\varphi)v_k(x^m|\varphi)\mu(dx^m).$$

From **[Mod]**, one can apply Theorem 4 (Nguyen and Zessin (1979)). It follows that there exists $n_0$ such that for all $n\ge n_0$

$$\left|\left(\underline{\boldsymbol{U}}_n^{(2)}(\varphi;\boldsymbol{\theta})\right)_{j,k}\right| \le 2\times e^K\boldsymbol{E}\left(\left|v_j(0^M|\Phi)v_k(0^M|\Phi)\right|\right) \text{ for } P_{\boldsymbol{\theta}^\star}-a.e.\ \varphi.$$

Note that from Theorem 4 (Nguyen and Zessin (1979)), for all $\boldsymbol{\theta}$ (and in particular $\boldsymbol{\theta}=\boldsymbol{\theta}^\star$), $\underline{\boldsymbol{U}}_n^{(2)}(\cdot;\boldsymbol{\theta})$ converges almost surely as $n\to+\infty$ towards $\underline{\boldsymbol{U}}^{(2)}(\boldsymbol{\theta})$ which is the $(p,p)$ matrix with elements

$$\left(\underline{\boldsymbol{U}}^{(2)}(\boldsymbol{\theta})\right)_{j,k} = \boldsymbol{E}\left(v_j(0^m|\Phi)v_k(0^m|\Phi)e^{-V(0^m|\varphi;\boldsymbol{\theta})}\right),\ \text{for } j,k=1,\dots,p.$$

Let us underline that $\underline{\boldsymbol{U}}^{(2)}(\boldsymbol{\theta})$ is a symmetric definite positive matrix. Indeed, it is a positive matrix since for all $\boldsymbol{y}\in\mathbb{R}^{p+1}$

$$\begin{aligned} \boldsymbol{y}^T\underline{\boldsymbol{U}}^{(2)}(\boldsymbol{\theta})\boldsymbol{y} &= \sum_{j,k} y_j\boldsymbol{E}\left(v_j(0^M|\Phi)v_k(0^M|\Phi)e^{-V\left(0^M|\Phi;\boldsymbol{\theta}\right)}\right)y_k \\ &= \boldsymbol{E}\left(\left(\boldsymbol{y}^T\boldsymbol{v}\left(0^M|\Phi\right)\right)^2 e^{-V\left(0^M|\Phi;\boldsymbol{\theta}\right)}\right) \ge 0. \end{aligned}$$

And it is definite since, for all $\boldsymbol{y}\in\mathbb{R}^{p+1}\setminus\{0\}$ from **[Ident]**, $\boldsymbol{y}^T\boldsymbol{v}\left(0^m|\varphi\right)=0$ for $\lambda^{\mathrm{m}}\otimes P_{\boldsymbol{\theta}^\star}-$a.e. $(m,\varphi)$ implies $\boldsymbol{y}=0$.

*Conclusion Step* Under the assumptions **[Mod]** and **[Ident]**, and using Steps 1 and 2, one can apply a classical result concerning asymptotic normality for minimum contrast estimators, see Guyon (1995) in order to obtain (4.3).

Now, the result (4.5) is proved in three substeps:

($i$) We first prove that the matrix $\underline{\boldsymbol{\Sigma}}(\boldsymbol{\theta}^\star) = \underline{\boldsymbol{\Sigma}}(\overline{D}, \boldsymbol{\theta}^\star)$ is a symmetric definite positive matrix. From Equation (6.9), it is sufficient to prove that the matrix $\mathbb{V}\mathrm{ar}(|\Lambda_n(\overline{D})|^{-1/2}\boldsymbol{LPL}^{(1)}_{\Lambda_n(\overline{D})}(\Phi;\boldsymbol{\theta}^\star))$ is definite positive for $n$ large enough. Let $\boldsymbol{y} \in \mathbb{R}^p \setminus \{0\}$, the aim is to prove that

$$V := \boldsymbol{y}^T \mathbb{V}\mathrm{ar}\left(|\Lambda_n(\overline{D})|^{-1/2}\boldsymbol{LPL}^{(1)}_{\Lambda_n(\overline{D})}(\Phi;\boldsymbol{\theta}^\star)\right)\boldsymbol{y} > 0.$$

Let $\overline{\Lambda} = \cup_{i\in\mathbb{B}\left(0,\left\lceil\frac{D}{\overline{D}}\right\rceil\right)}\Delta_i(\overline{D})$, using the same argument of Jensen and Künsch (1994) (Equation (3.2)), one can write

$$V \geq |\Lambda_n(\overline{D})|^{-1}\,\boldsymbol{E}\left(\mathbb{V}\mathrm{ar}\left(\boldsymbol{y}^T\boldsymbol{LPL}^{(1)}_{\Lambda_n(\overline{D})}(\Phi;\boldsymbol{\theta}^\star)\,|\Phi_{\Delta_k(\overline{D})}, k\notin(2\left\lceil\frac{D}{\overline{D}}\right\rceil+1)\mathbb{Z}^2\right)\right).$$

Note that for $i \neq j \in I_n$,

$$\begin{aligned}&Cov\left(\boldsymbol{y}^T\boldsymbol{LPL}^{(1)}_{\Delta_i(\overline{D})}(\varphi;\boldsymbol{\theta}^\star), \boldsymbol{y}^T\boldsymbol{LPL}^{(1)}_{\Delta_j(\overline{D})}(\varphi;\boldsymbol{\theta}^\star)\,|\Phi_{\Delta_k(\overline{D})}, k\notin(2\left\lceil\frac{D}{\overline{D}}\right\rceil+1)\mathbb{Z}^2\right)\\&\quad= 0\end{aligned}$$

due to the independence of $\boldsymbol{LPL}^{(1)}_{\Delta_i(\overline{D})}(\varphi;\boldsymbol{\theta}^\star)$ and $\boldsymbol{LPL}^{(1)}_{\Delta_j(\overline{D})}(\varphi;\boldsymbol{\theta}^\star)$ conditionally on $\Phi_{\Delta_k(\overline{D})}, k \notin (2\left\lceil\frac{D}{\overline{D}}\right\rceil+1)\mathbb{Z}^2$ when $i,j \in I_n \cap (2\left\lceil\frac{D}{\overline{D}}\right\rceil+1)\mathbb{Z}^2$ and $\boldsymbol{LPL}^{(1)}_{\Delta_i(\overline{D})}(\varphi;\boldsymbol{\theta}^\star)$ or $\boldsymbol{LPL}^{(1)}_{\Delta_j(\overline{D})}(\varphi;\boldsymbol{\theta}^\star)$ is constant when either $i$ or $j \notin I_n \cap (2\left\lceil\frac{D}{\overline{D}}\right\rceil+1)\mathbb{Z}^2$. As a direct consequence,

$$\begin{aligned}V &\geq |\Lambda_n(\overline{D})|^{-1}\boldsymbol{E}\left(\mathbb{V}\mathrm{ar}\left(\boldsymbol{y}^T\sum_{i\in I_n}\boldsymbol{LPL}^{(1)}_{\Delta_i(\overline{D})}(\Phi;\boldsymbol{\theta}^\star)\,\middle|\Phi_{\Delta_k(\overline{D})}, k\notin(2\left\lceil\frac{D}{\overline{D}}\right\rceil+1)\mathbb{Z}^2\right)\right)\\
&= |\Lambda_n(\overline{D})|^{-1}\sum_{i\in I_n}\boldsymbol{E}\left(\mathbb{V}\mathrm{ar}\left(\boldsymbol{y}^T\boldsymbol{LPL}^{(1)}_{\Delta_i(\overline{D})}(\Phi;\boldsymbol{\theta}^\star)\,\middle|\Phi_{\Delta_k(\overline{D})}, k\notin(2\left\lceil\frac{D}{\overline{D}}\right\rceil+1)\mathbb{Z}^2\right)\right)\\
&= |\Lambda_n(\overline{D})|^{-1}\sum_{\ell\in I_n\cap(2\left\lceil\frac{D}{\overline{D}}\right\rceil+1)\mathbb{Z}^2\setminus\widetilde{I_n}}\\
&\quad\boldsymbol{E}\left(\mathbb{V}\mathrm{ar}\left(\boldsymbol{y}^T\sum_{i\in I_n\cap\mathbb{B}\left(\ell,\left\lceil\frac{D}{\overline{D}}\right\rceil\right)}\boldsymbol{LPL}^{(1)}_{\Delta_i(\overline{D})}(\Phi;\boldsymbol{\theta}^\star)\,\middle|\Phi_{\Delta_k(\overline{D})}, k\notin(2\left\lceil\frac{D}{\overline{D}}\right\rceil+1)\mathbb{Z}^2\right)\right)\\
&\quad+|\Lambda_n(\overline{D})|^{-1}\sum_{\ell\in(2\left\lceil\frac{D}{\overline{D}}\right\rceil+1)\mathbb{Z}^2\cap\widetilde{I_n}}\\
&\quad\boldsymbol{E}\left(\mathbb{V}\mathrm{ar}\left(\boldsymbol{y}^T\sum_{i\in I_n\cap\mathbb{B}\left(\ell,\left\lceil\frac{D}{\overline{D}}\right\rceil\right)}\boldsymbol{LPL}^{(1)}_{\Delta_i(\overline{D})}(\Phi;\boldsymbol{\theta}^\star)\,\middle|\Phi_{\Delta_k(\overline{D})}, k\notin(2\left\lceil\frac{D}{\overline{D}}\right\rceil+1)\mathbb{Z}^2\right)\right)\end{aligned}$$

Following the proof of Step 1, condition $(iii)$ one may prove that the second right-hand term tends to 0 as $n \to +\infty$. Therefore by using the stationarity, we have for $n$ large enough

$$\begin{aligned} V &\geq \frac{1}{2}|\Lambda_n(\overline{D})|^{-1}\left|I_n \cap (2\left\lceil\frac{D}{\overline{D}}\right\rceil+1)\mathbb{Z}^2\right| \\ &\qquad \times \boldsymbol{E}\left(\mathbb{V}\text{ar}\left(\boldsymbol{y}^T\boldsymbol{LPL}_{\Lambda}^{(1)}(\Phi;\boldsymbol{\theta}^\star)\,\middle|\Phi_{\Delta_k(\overline{D})}, 1\leq|k|\leq 2\left\lceil\frac{D}{\overline{D}}\right\rceil\right)\right) \\ &= \frac{\overline{D}^{-2}}{2}\frac{|I_n\cap(2\left\lceil\frac{D}{\overline{D}}\right\rceil+1)\mathbb{Z}^2|}{|I_n|} \\ &\qquad \times \boldsymbol{E}\left(\mathbb{V}\text{ar}\left(\boldsymbol{y}^T\boldsymbol{LPL}_{\Lambda}^{(1)}(\Phi;\boldsymbol{\theta}^\star)\,\middle|\Phi_{\Delta_k(\overline{D})}, 1\leq|k|\leq 2\left\lceil\frac{D}{\overline{D}}\right\rceil\right)\right) \\ &\geq \frac{\overline{D}^{-2}}{2}\left(\frac{3}{4\left\lceil\frac{D}{\overline{D}}\right\rceil+1}\right)^2 \\ &\qquad \times \boldsymbol{E}\left(\mathbb{V}\text{ar}\left(\boldsymbol{y}^T\boldsymbol{LPL}_{\Lambda}^{(1)}(\Phi;\boldsymbol{\theta}^\star)\,\middle|\Phi_{\Delta_k(\overline{D})}, 1\leq|k|\leq 2\left\lceil\frac{D}{\overline{D}}\right\rceil\right)\right) \end{aligned}$$

Assume there exists some positive constant $c$ such that $P_{\boldsymbol{\theta}^\star}-$a.s. $\boldsymbol{y}^T\boldsymbol{LPL}_{\Lambda}^{(1)}(\Phi;\boldsymbol{\theta}^\star) = c$ when the variables $\Phi_{\Delta_k(\overline{D})}, 1 \leq |k| \leq 2\left\lceil\frac{D}{\overline{D}}\right\rceil$ are (for example) fixed to $\emptyset$. By assuming **[SDP]** it follows that for any $\varphi_i \in A_i$ for $i = 0,\ldots,\ell$ (with $\ell \geq p$), $\boldsymbol{y}^T\left(\boldsymbol{LPL}_{\Lambda}^{(1)}(\varphi_i;\boldsymbol{\theta}^\star) - \boldsymbol{LPL}_{\Lambda}^{(1)}(\varphi_0;\boldsymbol{\theta}^\star)\right) = 0$. Since for all $(\varphi_0,\ldots,\varphi_\ell) \in A_0\times\ldots\times A_\ell$, the matrix with entries $\left(\boldsymbol{LPL}_{\Lambda}^{(1)}(\varphi_i;\boldsymbol{\theta}^\star)\right)_j - \left(\boldsymbol{LPL}_{\Lambda}^{(1)}(\varphi_0;\boldsymbol{\theta}^\star)\right)_j$ is assumed to be injective, this leads to $\boldsymbol{y} = 0$ and hence to some contradiction. Therefore, when the variables $\Phi_{\Delta_k(\overline{D})}$, $1 \leq |k| \leq 2\left\lceil\frac{D}{\overline{D}}\right\rceil$ are fixed to $\emptyset$, the variable $\boldsymbol{y}^T\boldsymbol{LPL}_{\Lambda}^{(1)}(\Phi;\boldsymbol{\theta}^\star)$ is almost surely not a constant. Hence, $\underline{\boldsymbol{\Sigma}}(\boldsymbol{\theta}^\star)$ is a symmetric definite positive matrix.

$(ii)$ *Convergence of* $\widehat{\underline{\boldsymbol{\Sigma}}}_n(\varphi;D^\vee,\widetilde{D},\boldsymbol{\theta})$.

Let us recall that for any $\varphi\in\Omega$, $D^\vee \geq D$ and $\boldsymbol{\theta}\in\boldsymbol{\Theta}$ we define

$$\widehat{\underline{\boldsymbol{\Sigma}}}_n(\varphi;D^\vee,\widetilde{D},\boldsymbol{\theta}) = \frac{\widetilde{D}^{-2}}{|I_n|}\sum_{i\in I_n}\sum_{j\in I_n\cap\mathbb{B}\left(i,\left\lceil\frac{D^\vee}{\widetilde{D}}\right\rceil\right)}\boldsymbol{LPL}_{\Delta_i}^{(1)}(\varphi;\boldsymbol{\theta})^T\boldsymbol{LPL}_{\Delta_j}^{(1)}(\varphi;\boldsymbol{\theta})$$

We also define

$$X_i(\varphi) := X_i(\varphi)^{k,\ell} = \sum_{j\in I_n\cap\mathbb{B}\left(i,\left\lceil\frac{D^\vee}{\widetilde{D}}\right\rceil\right)}\left(\boldsymbol{LPL}_{\Delta_i}^{(1)}(\varphi;\boldsymbol{\theta})\right)_k\left(\boldsymbol{LPL}_{\Delta_j}^{(1)}(\varphi;\boldsymbol{\theta})\right)_\ell,$$

$Y_i(\varphi) := X_i(\varphi) - \boldsymbol{E}(X_i(\Phi))$ and $\overline{Y}_n(\varphi) = |I_n|^{-1}\sum_{i\in I_n}Y_i(\varphi)$. Since one may notice that $\boldsymbol{E}(X_i(\Phi)) = \widetilde{D}^2\left(\underline{\boldsymbol{\Sigma}}(\widetilde{D},\boldsymbol{\theta})\right)_{k,\ell}$, we have

$$\overline{Y}_n(\varphi) = \widetilde{D}^2\left(\widehat{\underline{\boldsymbol{\Sigma}}}_n(\varphi;D^\vee,\widetilde{D},\boldsymbol{\theta}) - \underline{\boldsymbol{\Sigma}}(\widetilde{D},\boldsymbol{\theta})\right)_{k,\ell}.$$

Thus, the aim is to prove that as $n \to +\infty$, $\overline{Y}_n(\varphi) \to 0$ for $P_{\boldsymbol{\theta}^\star}$−a.e. $\varphi$. Since the process $\{Y_i, i \in \mathbb{Z}^2\}$ is stationary, it is sufficient to prove, see *e.g.* Guyon (1995)

- $(a)$ $\boldsymbol{E}\left(Y_0(\Phi)^2\right) < +\infty$
- $(b)$ $\boldsymbol{E}\left(|I_n|\overline{Y}_n(\Phi)^2\right) < +\infty$.

$(a)$ We leave the reader to verify that **[Mod]** ensures this integrability condition.

$(b)$ Note that $Y_i(\varphi)$ depends only on $\varphi_{\Delta_j}$ for $j \in \mathbb{B}\left(i, \left\lceil \frac{D^\vee}{\widetilde{D}} \right\rceil + \left\lceil \frac{D}{\widetilde{D}} \right\rceil\right)$. Hence, by choosing $j \in I_n \cap \mathbb{B}\left(i, \alpha_{D^\vee,D}\right)^c$ with $\alpha_{D^\vee,D} = \alpha_{D^\vee,D}(\widetilde{D}) := 2\left\lceil \frac{D^\vee}{\widetilde{D}} \right\rceil + \left\lceil \frac{D}{\widetilde{D}} \right\rceil$, the variables $Y_i$ and $Y_j$ are independent. Then, we obtain

$$\begin{aligned}
\boldsymbol{E}\left(|I_n|\overline{Y}_n(\Phi)^2\right) &= \frac{1}{|I_n|}\sum_{i,j\in I_n} \boldsymbol{E}\left(Y_i(\Phi)Y_j(\Phi)\right) \\
&= \frac{1}{|I_n|}\sum_{i\in I_n}\left(\sum_{j\in I_n\cap\mathbb{B}\left(i,\alpha_{D^\vee,D}\right)} \boldsymbol{E}\left(Y_i(\Phi)Y_j(\Phi)\right)\right. \\
&\qquad \left. + \sum_{j\in I_n\cap\mathbb{B}\left(i,\alpha_{D^\vee,D}\right)^c} \boldsymbol{E}\left(Y_i(\Phi)Y_j(\Phi)\right)\right) \\
&= \frac{1}{|I_n|}\sum_{i\in I_n}\sum_{j\in I_n\cap\mathbb{B}\left(i,\alpha_{D^\vee,D}\right)} \boldsymbol{E}\left(Y_i(\Phi)Y_j(\Phi)\right) \\
&\sim \sum_{k\in\mathbb{B}\left(0,\alpha_{D^\vee,D}\right)} \boldsymbol{E}\left(Y_0(\Phi)Y_k(\Phi)\right) \leq (2\alpha_{D^\vee,D}+1)\,\boldsymbol{E}\left(Y_0(\Phi)^2\right).
\end{aligned}$$

Therefore, for all $D^\vee \geq D$ and for all $\boldsymbol{\theta} \in \boldsymbol{\Theta}$, we have for $P_{\boldsymbol{\theta}^\star}$−a.e. $\varphi$ as $n \to +\infty$

$$\widehat{\underline{\boldsymbol{\Sigma}}}_n(\varphi; D^\vee, \widetilde{D}, \boldsymbol{\theta}) \to \underline{\boldsymbol{\Sigma}}(\widetilde{D}, \boldsymbol{\theta}) = \underline{\boldsymbol{\Sigma}}(\boldsymbol{\theta}). \tag{6.10}$$

$(iii)$ Since for any $\varphi$, the functions $\underline{\boldsymbol{U}}_n^{(2)}(\varphi; \cdot)$ and $\widehat{\underline{\boldsymbol{\Sigma}}}_n(\varphi; D^\vee, \widetilde{D}, \cdot)$ are continuous, it follows from Step 2 and (6.10) that one obtains for $P_{\boldsymbol{\theta}^\star}$−a.e. $\varphi$, as $n \to +\infty$

$$\underline{\boldsymbol{U}}_n^{(2)}(\varphi; \widehat{\boldsymbol{\theta}}) \to \underline{\boldsymbol{U}}^{(2)}(\boldsymbol{\theta}^\star) \quad \text{and} \quad \widehat{\underline{\boldsymbol{\Sigma}}}_n(\varphi; D^\vee, \widetilde{D}, \widehat{\boldsymbol{\theta}}) \to \underline{\boldsymbol{\Sigma}}(\boldsymbol{\theta}^\star).$$

Finally, note that the previous convergence also implies that for $n$ large enough $\widehat{\underline{\boldsymbol{\Sigma}}}_n(\Phi; D^\vee, \widetilde{D}, \widehat{\boldsymbol{\theta}})$ is almost surely a symmetric definite positive matrix.

## Acknowledgements

We are grateful to the referees and the associate editor for their helpful comments. We also thank Laurence Pierret for the time spent on improving the writing of this paper in English.